# Framework for Modeling and Optimization of On-Orbit Servicing Operations under Demand Uncertainties[1]

Tristan Sarton du Jonchay[2], Hao Chen[3], Onalli Gunasekara[4]

and Koki Ho[5]

*Georgia Institute of Technology, Atlanta, GA, 30332, USA*

**This paper develops a framework that models and optimizes the operations of complex on-orbit servicing infrastructures involving one or more servicers and orbital depots to provide multiple types of services to a fleet of geostationary satellites. The proposed method extends the state-of-the-art space logistics technique by addressing the unique challenges in on-orbit servicing applications and integrates it with the Rolling Horizon decision making approach. The space logistics technique enables modeling of the on-orbit servicing logistical operations as a Mixed-Integer Linear Program whose optimal solutions can efficiently be found. The Rolling Horizon approach enables the assessment of the long-term value of an on-orbit servicing infrastructure by accounting for the uncertain service needs that arise over time among the geostationary satellites. Two case studies successfully demonstrate the effectiveness of the framework for (1) short-term operational scheduling and (2) long-term strategic decision making for on-orbit servicing architectures under diverse market conditions.**

## Nomenclature

$b_{svt}$ = Servicer dispatch variables

---

[2] Ph.D. Student, Daniel Guggenheim School of Aerospace Engineering, Atlanta, GA, AIAA Student Member.

[3] Ph.D. Student, Daniel Guggenheim School of Aerospace Engineering, Atlanta, GA, AIAA Student Member.

[4] Ph.D. Student, Daniel Guggenheim School of Aerospace Engineering, Atlanta, GA, AIAA Student Member.

[5] Assistant Professor, Daniel Guggenheim School of Aerospace Engineering, Atlanta, GA, AIAA Senior Member.

$c_s^{\text{delay}}$ = Penalty fees per unit time for providing a service with delays

$c_v^{\text{depot}}$ = Costs per unit time for operating an orbital depot

$c^{\text{launch}}$ = Launch cost per unit mass

$c_k^{\text{pdm}}$ = Purchase, development and manufacturing costs per unit mass of the non-vehicle commodity $k$

$c_v^{\text{pdm}}$ = Purchase, development and manufacturing costs per unit of the vehicle commodity $v$

$c_v^{\text{serv}}$ = Costs per unit time for operating a servicer

$\boldsymbol{c}_{\boldsymbol{vijt}}^{+}$ = Cost parameters for the non-vehicle commodities

$c'^{+}_{\boldsymbol{vijt}}$ = Cost parameters for the vehicle commodity

$\boldsymbol{d}_{it}$ = Demand parameters for the non-vehicle commodities

$d'_{ivt}$ = Demand parameters for the vehicle commodity

$\boldsymbol{e}_v$ = Design parameters of vehicle *v*

$H_{vij}$ = Concurrency constraint matrix

$h_{sv\tau}$ = Service assignment variables

$i, j$ = Node indices

$J$ = Objective function (*e.g.*, profits)

$J_{\text{delay}}$ = Component of the objective function capturing the penalty fees due to service delays

$J_{\text{depots}}$ = Component of the objective function capturing the operating costs of the depots

$J_{\text{launch}}$ = Component of the objective function capturing the launch costs

$J_{\text{pdm}}$ = Component of the objective function capturing the purchase, development, and manufacturing costs

$J_{\text{revenues}}$ = Component of the objective function capturing the revenues generated by the services

$J_{\text{serv}}$ = Component of the objective function capturing the operating costs of the servicers

$\mathcal{K}_C$ = Index set for continuous commodities

$\mathcal{K}_I$ = Index set for integer commodities

$k$ = Non-vehicle commodity index

$\mathcal{K}_{\text{tools}}$ = Index set for the servicers' tools

$\mathcal{L}_{ij}$ = Index set for the time steps associated with mission operation opportunities

$m_k$ = Unit mass for the non-vehicle integer commodity $k$

$m_v$ = Unit mass for the vehicle commodity $v$

$\mathcal{N}$ = Index set for all nodes in the static network

$\mathcal{N}_c$ = Index set for the customer nodes in the static network

$\mathcal{N}_e$ = Index set for the Earth nodes in the static network (*i.e.*, spaceports)

$\mathcal{N}_p$ = Index set for the on-orbit servicing parking nodes in the static network

$n$ = Number of spaceflight time steps per time interval $T$ in the time-expanded network

$Q_{vij}$ = Commodity transformation matrix

$r_s$ = Revenues received upon providing service $s$

$\mathcal{S}$ = Index set for the service needs

$\mathcal{S}_i$ = Index set for the service needs arising at customer node $i$

$s$ = Service index

$s_v$ = Structural mass of vehicle $v$

$\mathcal{T}$ = Index set for the time steps

$T$ = Service time interval in time-expanded network

$t, \tau$ = Time indices

$\Delta t$ = Length of the spaceflight time steps in the time-expanded network

$\Delta_{ijt}$ = Duration of the arc between node $i$ and node $j$ at time $t$

$\mathcal{V}$ = Index set for the vehicles

$\mathcal{V}_{\text{depot}}$ = Index set for the orbital depots

$\mathcal{V}_{\text{serv}}$ = Index set for the servicers

$v$ = Vehicle index

$\mathcal{W}_s$ = Index set for the time steps within the service window associated with service $s$

$\boldsymbol{x}_{vijt}^{\pm}$ = Vector of commodity flow variables for the non-vehicle commodities

$x_{vijtk}^{\pm}$ = Individual commodity flow variables for the non-vehicle commodity $k$

$y_{vijt}^{\pm}$ = Commodity flow variables for the vehicle commodity $v$

$\beta_{s\tau t}$ = Binary parameters relating servicer dispatch variables and service assignment variables

$\gamma_{sk}$ = Service-tool mapping between service $s$ and tool $k$

$\tau_s$ = Date of occurrence of service need $s$

## I. Introduction

Geosynchronous-equatorial-orbit (GEO) satellites have played a critical role as a key infrastructure for our life, but they are also known for the difficulty of their maintenance due to their high altitude. Recently, a decrease in the replacement of GEO satellites has been observed [1]. This indicates that, while GEO satellite operators are weighing the pros and cons between investing in large-scale low-Earth-orbit (LEO) constellations and investing in high-throughput GEO satellite technologies [2], the fleet of operational GEO satellites is aging. This is a unique opportunity for On-Orbit Servicing (OOS) providers to step in and sell life extension services to GEO satellite operators at the fraction of a cost of a replacement satellite. In addition, recent trends in semi-autonomous and autonomous robotic research (*e.g.*, DARPA's Robotic Servicing of Geosynchronous Satellites (RSGS) [3]) show that reliable OOS operations are becoming a reality, making satellite operators more confident in satellite servicing. Intelsat for instance has contracted Northrop Grumman to have their first Mission Extension Vehicle (MEV-1) extend the life of Intelsat-901 for five years [4].

Despite the favorable turn that OOS is taking, this industry is still very much at its infancy. The tradespace of OOS is very large thus its exploration has been a challenging task. The existing GEO servicing literature can first be classified based on the type of OOS mission concept investigated:

(1) *One-to-one* architectures in which one servicer (the servicing spacecraft) is serving only one client satellite [5]. The first OOS missions will likely involve one-to-one architectures with servicers designed to provide one type of service such as station keeping. Northrop Grumman's on-going MEV-1 mission [4] is a one-to-one example. Infrastructure elements such as orbital depots may also be considered to further support the OOS operations;

(2) *One-to-many* architectures in which one servicer provides services to a fleet of several client satellites [6-11]. An orbital spare depot or refueling station may be considered to further support the OOS operations;

(3) *Many-to-many* architectures in which several servicers provide services to a fleet of client satellites [12-15]. Multiple orbital spare depots and refueling stations may be considered to further support the OOS operations.

For each OOS mission concept, the literature may then further be classified in two categories:

a. The analysis and design of accurate orbital trajectories associated with multi-transfer servicing operations, using either high-thrust [6] or low-thrust servicers [12]. These analyses usually consider high-fidelity models of perturbation forces to design the transfer trajectories;
b. The analyses of the operational strategy of a particular OOS infrastructure from a scheduling [7] and/or design standpoint [5,8,9]. Simulation [8,10,13], optimization [7,9,11] or a mixture of both [5] may be used to explore the value of the investigated OOS architecture by varying its operational scheduling and/or design.

Despite a prolific OOS literature, the operational scheduling and system design of sustainable many-to-many OOS infrastructures has been largely unexplored. This is because this problem is highly complex and requires a sophisticated optimization problem formulation to handle that.

Leveraging the recent development in space logistics, this paper aims to fill that gap by developing a new framework to optimize the OOS operations and system design under service demand uncertainties. The developed method will extend the state-of-the-art Mixed-Integer Linear Programming (MILP)-based space logistics modeling to address the unique challenges in OOS applications and integrate it with the dynamic Rolling Horizon (RH) decision making approach. The resulting framework enables both the short-term scheduling of OOS operations and the long-term decision making for the OOS systems design.

Space logistics have largely been explored in the literature to optimize the design and planning of complex human space exploration missions through heuristics methods [16], simulations [17], graph theory [18], or network flow models [19-24]. The method described in this paper extends the latter to accurately model OOS operations. In the traditional space logistics network flow models, the formulation only specifies as inputs the demand of commodities at locations of interest such as Mars or the Moon. The optimizer then decides when and where to assign space vehicles in order to satisfy the demand. In the context of OOS, the servicers must stay a certain amount of time at the satellites' locations to provide the services. The traditional space logistics formulation does not capture this constraint and therefore a unique extension is needed for the OOS applications. The advantage of network flow models is that they can be modeled as MILP problems which can theoretically reach global optima.

RH decision making, which is leveraged to design a long-term strategy for OOS, has been typically used to make decisions in a dynamic stochastic environment, usually characterized by uncertainties in demand and the resulting cost in forecasting this demand. More specifically, this technique consists in making a series of decisions repeatedly based on a short- to mid-term forecast of future demand [25,26]. An advantage of RH is to decompose a complex dynamic scheduling problem into smaller sub-problems whose combined optimal solutions yield a satisfying solution of the original problem at a lesser computational cost. In the literature, RH has been used for the optimal scheduling of Earth observing satellites [27] but never for the operational scheduling of OOS infrastructures.

The OOS mission concept explored in this paper is as follows. A servicing company launches and operates orbital depots and impulsive-thrust servicers to provide services to a fleet of GEO satellites whose service needs occur on a random basis (*e.g.*, component failure) or on a deterministic basis (*i.e.*, pre-planned such as propellant depletion). Whenever a satellite needs a service, the servicing company first decides whether to provide the service, and if so, a servicer is dispatched to the client satellite to start the service before a certain deadline. Once the service is done, the servicer is available for dispatch to a new service. The notional services considered in this paper are inspection, refueling, station keeping, satellite repositioning, repair, mechanism deployment, and retirement. The servicers and orbital depots may be replenished on a regular basis with launch vehicles sent from Earth.

In order to make OOS ventures profitable and competitive, the scheduling of their operations as well as the associated supply chain of commodities need to be optimized concurrently. The optimization framework proposed in this paper is a flexible tool that OOS companies could leverage to devise the best short-term and long-term operational strategies given the uncertain service demand. This is enabled through three novelties summarized here: (1) modeling the logistics of complex many-to-many OOS infrastructures as a network; (2) extending the traditional space logistic MILP formulation with additional variables and constraints to make the OOS operations realistic; (3) leveraging RH decision making to explore the long-term strategic planning of sustainable OOS infrastructures.

The rest of this paper is organized as follows. Section II presents the modeling of the customer satellites, the OOS infrastructure, and the OOS logistics network. Section III describes the developed methods based on the space logistics and RH approaches to optimize the OOS operations. Section IV demonstrates the value of the proposed methods for applications such as the short-term operational scheduling and long-term strategic planning of complex many-to-many OOS infrastructures. Finally, section V concludes the paper.

## II. Modeling and Network Overview

This section presents the general modeling of the fleet of customer satellites and their service demands (subsection II.A), the modeling of the OOS architecture (subsection II.B), and an overview of the OOS logistics static network (subsection II.C).

### A. Customer Satellites and Orbital Service Needs

In this paper, the customer satellites are assumed to be distributed along the same circular orbit. The orbital states of the customer satellites thus only differ by their positions on the shared orbit (*e.g.,* the longitudes for GEO satellites). The methods developed in this paper apply to GEO servicing, which is one of the largest potential markets for OOS, as well as to the servicing of a fleet of customer satellites on any other circular orbit.

Other parameters of the customer satellites include their masses, the type of propellant they use for station keeping, and what type of life extension service they require (*e.g.,* refueling or station keeping) if such a need arises during their lifetimes. In addition, each defined satellite is specified parameters that model the characteristics of their potential service needs.

Seven types of service needs are considered in this paper. Three of them are assumed deterministic, which occur on a regular basis, whereas the remaining four are assumed to arise randomly. Tables 1 and 2 define the deterministic and random service needs modeled in this work. The main parameters used to model the deterministic and random service needs are defined in Table 3. Section IV assigns to these parameters data that is mostly drawn from the literature.

**Table 1. Definition of the deterministic service needs (*i.e.,* that can be planned ahead of time).**

| | Inspection | Refueling | Station Keeping |
|---|---|---|---|
| **Definition** | The servicer performs a proximity maneuver near the satellite, and inspects its condition without docking to it. | The servicer rendezvouses and docks to the satellite to top up its tank with additional propellant. | The servicer rendezvouses and docks to the satellite to perform station-keeping maneuvers in place of the satellite. |

**Table 2. Definition of the random service needs (*i.e.*, that cannot be planned ahead of time).**

| | Repositioning | Retirement | Repair | Mechanism Deployment |
|---|---|---|---|---|
| **Definition** | The servicer rendezvouses and docks to the satellite to change its orbital slot. | The servicer rendezvouses and docks to the defunct satellite to transport it to some graveyard orbit. | The servicer rendezvouses and docks to the satellite and replaces the defective parts with spare ones. | The servicer rendezvouses and docks to the satellite and unlocks its stuck appendages. |

**Table 3. Definitions of the main parameters used to model the service needs.**

| Parameter | Applies to deterministic and/or random needs? | Definition |
|---|---|---|
| *Revenue* | Deterministic and random | Revenue received by the OOS infrastructure for providing a service in response to some service need. |
| *Delay penalty cost* | Deterministic and random | How much it costs daily to the OOS operator to delay the beginning date of the service. |
| *Service duration* | Deterministic and random | The time it takes to provide a service. |
| *Service window* | Deterministic and random | Time interval within which the service must start, provided that the OOS operator decides to provide the service. |
| *Interoccurrence time* | Deterministic only | Used to generate service needs regularly spaced in time. |
| *Mean interoccurrence time* | Random only | Used to generate service needs according to a Poisson probability distribution. |

### B. On-Orbit Servicing Architectural Elements

The OOS infrastructures modeled in this paper include the servicers and their tools, orbital depots, and launch vehicles. A servicer is defined here as a fully- or semi-autonomous robotic spacecraft that uses tools to provide services to customer satellites. An OOS infrastructure needs at least one such servicer to physically access the customer satellites and provide them services. The model of a servicer is captured by four sets of parameters defined in Table 4. A fundamental assumption made in this paper regarding the servicers is that they use impulsive propulsion technologies. Thus, the rocket equation can be integrated in a MILP formulation as a linear equation relating the total servicer mass and its propellant consumption. Future developments of this OOS optimization framework will relax this assumption by modeling low-thrust servicers and their trajectories. The orbital maneuvers considered in this paper are phasing maneuvers under altitude and time-of-flight constraints. The way they are designed is presented in Appendix.

The servicers must also have the right tool(s) to provide the services. A mapping between the service needs presented in subsection II.A and four notional tools is presented in Table 5. Note that the optimization framework is flexible enough to consider any kind of user-defined tools. In this paper, we assume that a service cannot be provided by multiple types of tools, but some tools can be used to provide multiple types of services. The parameters used to model a tool are its mass and its cost (development, manufacturing, etc.).

**Table 4. Parameters used to model the servicers.**

| Parameter set | Definition | Examples of parameters |
|---|---|---|
| *Tool parameters* | Specify how many tools the servicer is equipped with and of which type those tools are. | Tools integrated to the servicer's bus (*cf.* Table 5) |
| *Orbital transfer parameters* | Specify the performance of the servicer in executing phasing maneuvers. | Specific impulse (*Isp*); propellant tank capacity; dry mass |
| *Cost parameters* | Specify how much it costs the OOS company to develop, manufacture, and operate a servicer | Development costs; operating costs |
| *Payload capacity parameters* | Capacities of the servicer's payload bays to provide services that require commodities (*e.g.*, refueling, repair) | Capacity of the servicer's tank used to transport customer satellites' propellant; payload capacity to transport spare parts |

**Table 5. Notional service-tool mapping (represented by an 'X' in the table).**

| | **T1: Refueling apparatus** | **T2: Observation sensors** | **T3: Dexterous robotic arm** | **T4: Capture mechanism** |
|---|---|---|---|---|
| *Inspection* | | **X** | | |
| *Refueling* | **X** | | | |
| *Station Keeping* | | | | **X** |
| *Repositioning* | | | | **X** |
| *Retirement* | | | | **X** |
| *Repair* | | | **X** | |
| *Mechanism Deployment* | | | **X** | |

In addition to the servicers and their tools, an OOS infrastructure may make use of one or more orbital depots that act as in-space propellant depots, spare warehouses, and robotic handling systems (*e.g.,* the International Space Station). Their functions are to store commodities needed to provide services, such as refueling and repair, and to transfer commodities between any two vehicles (*e.g.*, between the payload bay of a launch vehicle and a servicer). A depot offers the advantage of more responsive services but at the cost of a more expensive OOS infrastructure. The main parameters of a depot are its dry mass (to calculate how much it costs to deploy it in orbit), the capacities of its payload bays (*e.g.,* how much propellant and spares it can store), its manufacturing and operating costs, and its own propellant consumption in order to perform station keeping. Two important assumptions related to depots are made in this paper: (1) depots never execute orbital maneuvers, and (2) they are deployed at some orbital slots (referred to as *OOS parking nodes* in this paper) located on the same orbital ring as the fleet of customer satellites. This latter assumption facilitates the modeling of the OOS scheduling problem as a MILP as was already discussed in subsection II.A. Future developments of this framework will relax this assumption.

A launch vehicle is modeled with three main parameters: its launch frequency (assumed deterministic in this paper); its launch price tag (in $/kg); and its payload capacity to launch commodities to Geosynchronous Transfer Orbit (GTO) or GEO. In this paper, launch vehicles are assumed to rendezvous with depots and/or servicers at the OOS parking nodes.

Finally, the consumables considered in this paper are the spare parts, bipropellant, and monopropellant. Only one type of spare commodity is considered for simplicity in later analyses; this assumption can be relaxed. Multiple types of propellant were defined to allow more flexibility in the modeling of the servicers and customer satellites. For example, nowadays, GEO satellites use either monopropellant thrusters or ion engines for station keeping, while orbital maneuvers such as apogee kicks (to transfer from GTO to GEO) or inclination change are done primarily using monopropellant or bipropellant. In this paper, servicers are assumed to use bipropellant although this can easily be changed in the OOS optimization framework.

### C. Static Network Overview

As a transition to section III, this subsection introduces the static OOS logistics network based on which the MILP formulation is developed, while at the same time giving a high-level overview of the information introduced in subsections II.A and II.B.

The static network used for the optimization is made of three types of nodes: (1) the customer nodes, where the customer satellites are positioned and the service needs arise; (2) the OOS parking nodes, where orbital depots may be staged and servicers can get supplies from the depot and/or the launch vehicles; (3) the Earth nodes, where supplies of commodities are generated and launched to the OOS parking nodes via launch vehicles.

Direct arcs are then defined to allow the flow of commodities between any two nodes. Multiple arcs may be defined between two nodes to represent different vehicle options. A notional network is illustrated in Fig. 1 (not all arcs are represented for readability). Note that vehicles are considered as commodities in space logistics network flow models.

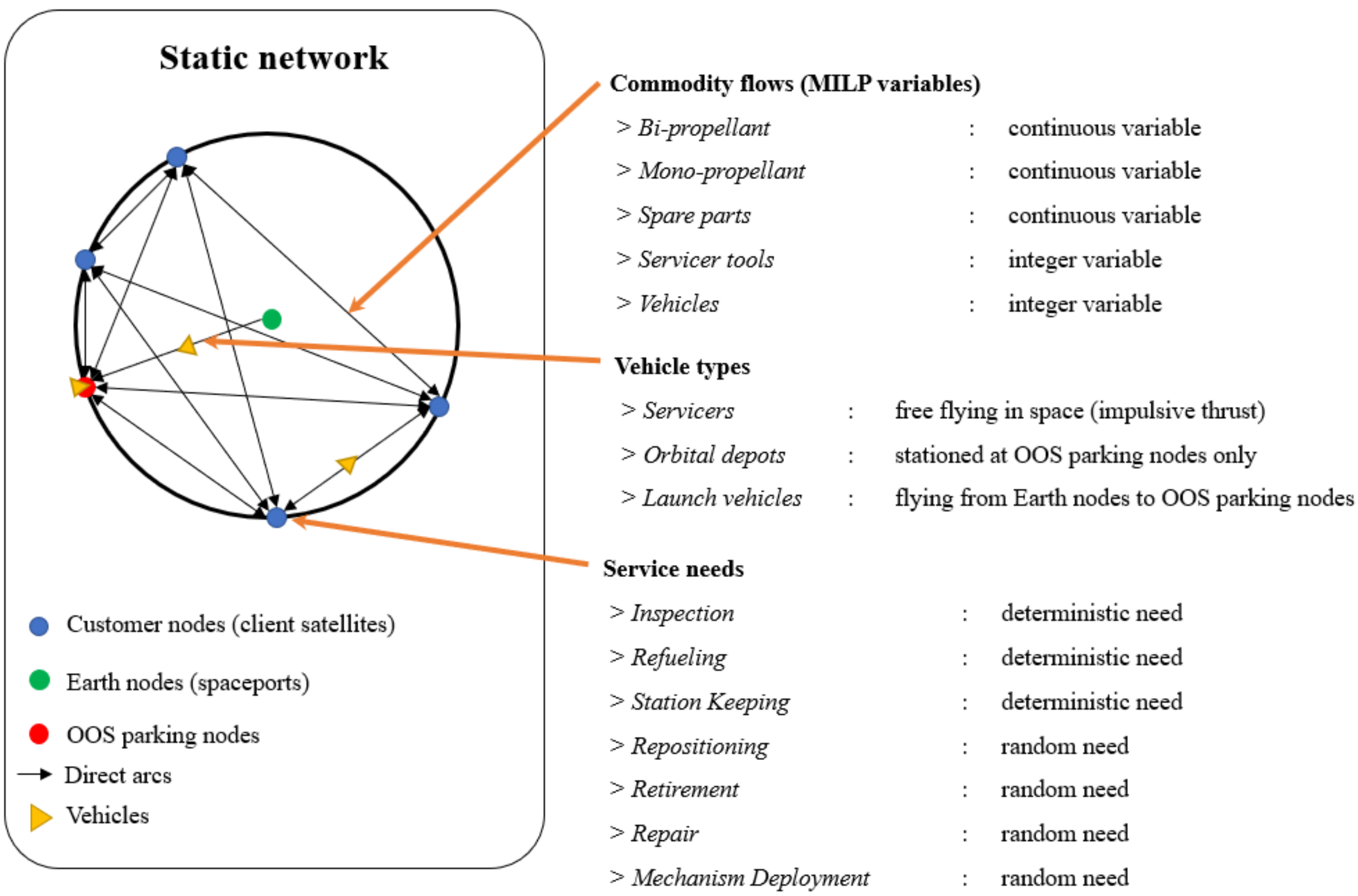

**Fig. 1 High-level overview of the static OOS logistics network used in developing the OOS optimization framework.**

## III. Methods

This section introduces the three important methodological pieces making up the proposed OOS optimization framework. It is desired to (1) accurately model the OOS operations over the time dimension, (2) formulate the OOS scheduling problem as a MILP of which a global optimum can theoretically be reached, given enough time to solve it, and (3) inform long-term decision making while considering uncertainties inherent to the random service needs (*e.g.,* repair). Subsection III.A describes how the static network presented in subsection II.C is expanded over time. Next, subsection III.B presents the traditional space logistics MILP formulation and its extension to accurately model OOS operations. Finally, Subsection III.C presents how the traditional RH decision making methodology is adapted in the context of OOS to allow for long-term sustainable architectural decisions under uncertainties.

### A. Static Network Time Expansion

This subsection shows how one can model the OOS operations over the time dimension by replicating the static OOS logistics network at predefined time steps. This may be done with a very fine time discretization at a cost of a large optimization time. Indeed, the number of MILP variables increases linearly with the number of time steps at which the static network is replicated. On the other hand, a very coarse time discretization may lead to inaccuracies in the OOS operations. Based on the application and the time available to run the analyses, a coarse or fine time discretization may be selected. This is done through the definition of three parameters: (1) the length $\Delta t$ of so-called *spaceflight time steps* that are used to accurately model the transportation of vehicles over the time-expanded network (also referred to as *dynamic network* in this paper); (2) a service time interval $T$; and (3) the number $n$ of spaceflight time steps within the time interval $T$. The following introduces the concept of each parameter.

The parameter $\Delta t$ controls the coarseness of the dynamic expansion of the static network. It is a divisor common to all time steps defined in the dynamic network.

The parameter $T$ is chosen at most as large as the greatest common divisor of the durations of the considered services. For example, if refueling services are provided in 20 days and inspection services in 30 days, then the greatest common divisor would be 10 days and we would need to choose $T \leq 10$ days. Note that the typical service durations found in the literature are 10-30 days [7,11] which allow for a naturally coarse time discretization. The static network is then replicated every time interval $T$.

The parameter $n$ defines the number of spaceflight time steps per interval $T$ that are used to accurately represent the finite flights of the servicers. They represent the first $\Delta t$-long time steps defined within each time interval $T$. The remaining time step of length $T - n \cdot \Delta t$, called *service time step*, is used exclusively for the provision of services, which means that vehicles cannot fly during that time step. Note that the impulsive trajectories are designed with a time of flight not exceeding $\Delta t$ justifying the considered time-expanded topology of the network.

Figure 2 depicts a notional time-expansion of a simple 3-node static network for $\Delta t = 1$ day, $T = 10$ days, and $n = 2$. This means that the servicers can fly over the network only between $t = 0$ and $t = 2$, between $t = 10$ and $t = 12$, etc. The time steps defined between $t = 2$ and $t = 10$, between $t = 12$ and $t = 20$, etc., are only used for the provision of services. This is illustrated by the bold, yellow direct arc on Fig. 2: between $t = 0$ and $t = 2$, a servicer is launched from Earth to the OOS parking node via a launch vehicle and then flies on its own from the OOS parking node to the customer node; between $t = 2$ and $t = 20$, the servicer provides a service before going back to the OOS parking node for storage at $t = 20$. Note that the service needs may only arise at $t = 0, 10, 20$, etc.

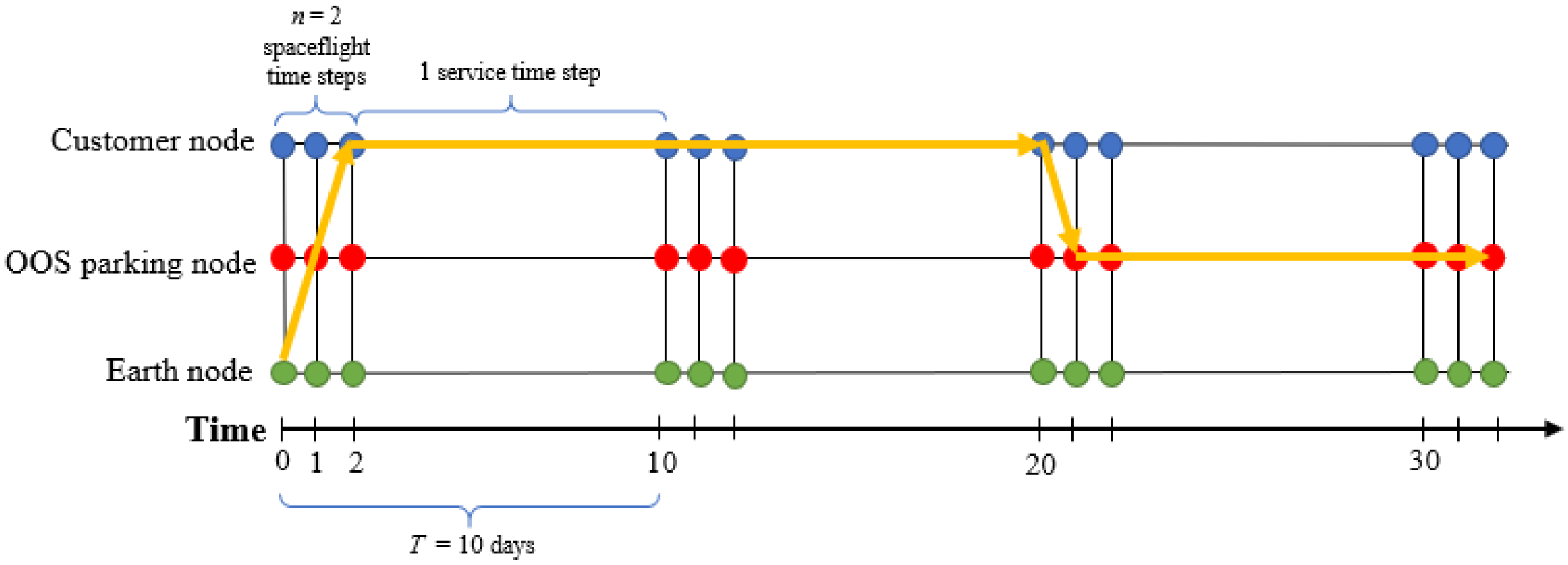

**Fig. 2 Notional time-expanded network and servicer path with $\Delta t$ = 1 day, $T$ = 10 days, and $n$ = 2.**

### B. Mathematical Formulation for On-Orbit Servicing Logistics

This subsection introduces the next methodological piece making up the proposed OOS optimization framework: the mathematical OOS logistics optimization formulation. In particular, we present the MILP model used for traditional space logistics problems (*e.g.*, interplanetary human space mission design), and then describe why that method cannot be used for OOS operations, followed by our developed extension for OOS application.

#### *1. Traditional Space Logistics Mixed Integer Linear Programming Formulation*

The traditional space logistics problem is formulated as a MILP model that builds upon the dynamic network formalism presented in Fig. 2. Let's consider a time-expanded network graph $\mathcal{G}$ made up of a set of nodes $\mathcal{N}$ and a set of direct arcs $\mathcal{A}$. The latter include both transportation arcs that connect different nodes and holdover arcs that connect a node with itself. Each arc has an index $(v, i, j)$, meaning that vehicle $v$ flies from node $i$ to node $j$. Additionally, the set of vehicles and time steps are denoted $\mathcal{V}$ and $\mathcal{T}$ respectively.

In the mathematical formulation presented here, the distinction is made between the vehicle commodity and the other commodities (*cf* Fig.1). The vehicle commodity flowing along arc $(v, i, j)$ at time step $t$ is denoted by the nonnegative scalar $y^{\pm}_{vijt}$. The definition is similar for the other commodities denoted by the $p$×1 column vector $\boldsymbol{x}^{\pm}_{vijt}$ of nonnegative scalars, where $p$ is the number of non-vehicle commodities. The commodities flow over the arcs and are split between node outflow $y^{+}_{vijt}/\boldsymbol{x}^{+}_{vijt}$ (*i.e.*, they leave node $i$) and node inflow $y^{-}_{vijt}/\boldsymbol{x}^{-}_{vijt}$ (*i.e.*, they enter node $j$). The vehicle commodity flow variables $y^{\pm}_{vijt}$ are integers. The non-vehicle commodity flow variables are either integer or continuous: the set $\mathcal{K}_I$ denotes the commodities that are integer, whereas $\mathcal{K}_C$ denotes the commodities that are continuous. Figure 1 shown earlier gives the natures of the variables that represent the commodities considered in this paper.

The cost coefficients are denoted by the scalar $c'^{+}_{vijt}$ for the vehicle commodity, and by the 1×$p$ row vector $\boldsymbol{c}^{+}_{vijt}$ for the other commodities. Each node $i$ of the static network is assigned at each time step a scalar vehicle demand $d'_{ivt}$ and a $p$×1 demand vector $\boldsymbol{d}_{it}$ for the other commodities, where demand and supply are assigned non-positive and non-negative values respectively. The MILP formulation presented below is based upon Ref. [22].

Minimize

$$J = \sum_{t \in \mathcal{T}} \sum_{(v,i,j) \in \mathcal{A}} \left\{ \boldsymbol{c}^{+}_{vijt} \boldsymbol{x}^{+}_{vijt} + c'^{+}_{vijt} y^{+}_{vijt} \right\} \tag{1}$$

Subject to:

$$\sum_{(v,j):(v,i,j) \in \mathcal{A}} \boldsymbol{x}^{+}_{vijt} - \sum_{(v,j):(v,j,i) \in \mathcal{A}} \boldsymbol{x}^{-}_{vji(t-\Delta_{jit})} \leq \boldsymbol{d}_{it} \qquad \forall t \in \mathcal{T}, \forall i \in \mathcal{N}, \forall v \in \mathcal{V} \tag{2}$$

$$\sum_{(v,j):(v,i,j) \in \mathcal{A}} y^{+}_{vijt} - \sum_{(v,j):(v,j,i) \in \mathcal{A}} y^{-}_{vji(t-\Delta_{jit})} \leq d'_{ivt} \qquad \forall t \in \mathcal{T}, \forall i \in \mathcal{N}, \forall v \in \mathcal{V} \tag{3}$$

$$Q_{vij} \begin{bmatrix} \boldsymbol{x}^{+}_{vijt} \\ s_v y^{+}_{vijt} \end{bmatrix} = \begin{bmatrix} \boldsymbol{x}^{-}_{vijt} \\ s_v y^{-}_{vijt} \end{bmatrix} \qquad \forall t \in \mathcal{T}, \forall (v,i,j) \in \mathcal{A} \tag{4}$$

$$H_{vij} \boldsymbol{x}^{+}_{vijt} \leq \boldsymbol{e}_v y^{+}_{vijt} \qquad \forall t \in \mathcal{T}, \forall (v,i,j) \in \mathcal{A} \tag{5}$$

$$\begin{cases} \boldsymbol{x}^{\pm}_{vijt} \geq 0_{p\times1} & \text{if } t \in \mathcal{L}_{ij} \\ \boldsymbol{x}^{\pm}_{vijt} = 0_{p\times1} & \text{otherwise} \end{cases} \qquad \forall t \in \mathcal{T}, \forall (v,i,j) \in \mathcal{A} \tag{6}$$

$$\boldsymbol{x}^{\pm}_{vijt} = \begin{bmatrix} x_1 \\ x_2 \\ \vdots \\ x_p \end{bmatrix}^{\pm}_{vijt}, \quad \begin{matrix} x_k \in \mathbb{R}_{\geq 0} \; \forall k \in \mathcal{K}_C \\ x_k \in \mathbb{Z}_{\geq 0} \; \forall k \in \mathcal{K}_I \end{matrix} \qquad \forall t \in \mathcal{T}, \forall (v,i,j) \in \mathcal{A} \tag{7}$$

$$y^{\pm}_{vijt} \in \mathbb{Z}_{\geq 0} \qquad \forall t \in \mathcal{T}, \forall (v,i,j) \in \mathcal{A} \tag{8}$$

Table 6 defines the indices, variables and parameters appearing in the MILP formulation. Table 7 explains the meaning of each equation appearing in the MILP formulation.

**Table 6. Definition of the indices, variables and parameters appearing in the space logistics MILP formulation.**

| | |
|---|---|
| **Indices** | |
| $v$ | Vehicle index. |
| $i, j$ | Node index. |
| $t$ | Time index. |
| **Variables** | |
| $\boldsymbol{x}_{vijt}^{\pm}$ | Non-vehicle commodity inflows/outflows. Commodities in $\boldsymbol{x}_{vijt}^{\pm}$ are integer or continuous variables based on the commodity type. For example, servicer tools are modeled as integers whereas propellant mass is modeled as continuous. |
| $y_{vijt}^{\pm}$ | Number of vehicles of type $v$ flying from node $i$ to node $j$ at time $t$. Integer variable. |
| **Parameters** | |
| $\boldsymbol{d}_{it}, d'_{ivt}$ | The demand of commodities and vehicles. In the OOS context, $\boldsymbol{d}_{it}$ represents the commodities needed to provide services (*e.g.*, propellant for refueling services). |
| $\Delta_{ijt}$ | The duration of the arc between node $i$ and node $j$ at time $t$. |
| $s_v$ | Structure mass of vehicle $v$. |
| $\boldsymbol{e}_v$ | Vehicle design parameters, including payload and propellant capacities. |
| $\boldsymbol{c}_{vijt}^{+}$ | Cost coefficients of the non-vehicle commodities. |
| $c'^{+}_{vijt}$ | Cost coefficients of the vehicle commodity. |
| $Q_{vij}$ | Commodity transformation matrix. |
| $H_{vij}$ | Concurrency constraint matrix. |

**Table 7. Meaning of the equations in the MILP formulation.**

| Equation | Name | Description |
|---|---|---|
| Objective (1) | Objective function | Minimize total mission cost. |
| Constraint (2) | Mass balance | The mass balance for the non-vehicle commodity flow into and out of node $i$. |
| Constraint (3) | Mass balance | The mass balance for the vehicle commodity flow into and out of node $i$. |
| Constraint (4) | Flow transformation | Transformation of the payloads during the flight along the arcs (*e.g.*, propellant burning). |
| Constraint (5) | Flow concurrency | Concurrency constraints of the flow related to the vehicles, such as payload bay and propellant tank capacities. |
| Constraint (6) | Flow bound | Time window (only when a time window opens is the commodity flow permitted). |

## 2. *Extension of the Traditional Space Logistics Formulation for On-Orbit Servicing Logistics*

Despite the effectiveness of the traditional space logistics formulation in optimizing complex interplanetary space missions [19-24], it has a caveat that makes it inapplicable to OOS operations. Namely, it cannot constrain a vehicle to remain at a node of the network. This constraint is not needed or is negligible for interplanetary space logistics mission design, but it will be critical when we consider the OOS operations because servicers must be providing services for some duration to the customer satellites.

To tackle this caveat, we extend the traditional space logistics formulation. This subsection first introduces two new sets of binary variables and the concept of *service window*. OOS-specific constraints are then defined to properly model OOS operations. Finally, the objective function used in this paper is presented.

### i. OOS-specific binary variables

The OOS-specific binary variables are introduced to allow the definition of additional constraints that properly assign services to servicers and force servicers to remain at customer satellites' locations wherever services must be provided. The two new sets of binary variables modeling this are the *service assignment variables* and the *servicer dispatch variables*.

#### *Service assignment variables and service window*

The service assignment variables, denoted $h_{sv\tau}$, are defined to specify at what date $\tau$ a given service $s$ must start and which servicer $v$ must provide it. More specifically, $h_{sv\tau} = 1$ if the optimizer requires servicer $v$ to start service $s$ at the date $\tau$. Before proceeding further with the description, multiple new index sets associated with the indices $v$, $s$, $\tau$ must be introduced. First, the set of servicers $\mathcal{V}_{\text{serv}}$, associated with index $v$, is a subset of the set of vehicles $\mathcal{V}$ introduced in Section III.B.1. Second, the set of services $\mathcal{S}$ includes all combined customer satellites' service needs occurring over the planning horizon. Also, the set of services $\mathcal{S}_i$ is a subset of $\mathcal{S}$ that corresponds to the services occurring at customer node $i$ over the planning horizon. (The concept of planning horizon will be defined more precisely in Section III.C.) Finally, the index set associated with $\tau$ is dependent upon the service index $s$ and is denoted $\mathcal{W}_s$. It comprises every time step defined in the dynamic network at which service $s$ may start. Thus, $\mathcal{W}_s$ is a subset of the set of time steps $\mathcal{T}$ introduced in Section III.B.1. It is related to a core concept of the proposed formulation, called *service window*, which is presented below.

Whenever a service need occurs at a given customer satellite, the OOS operator may be given some flexibility to select the date at which to dispatch a servicer to the satellite in need. This is allowed through the definition of a service window per each service need $s$ that arises over the planning horizon. For instance, let's assume that an OOS operator and a satellite operator agree on a 10-day service window whenever a refueling service need arises. Assume now that such a refueling service need, indexed by $s_0$, arises at that satellite at the date $t = 20$. The OOS operator decides to dispatch a servicer $v_0$ to provide the service $s_0$ at some date $\tau$ between $t = 20$ and $t = 30$. The optimizer embedded in the proposed framework would play the role of the OOS operator and decide at which date the service should start.

Figure 3 illustrates the above example with a simple dynamic network comprised only of the customer node where the service need $s_0$ arises. The $\Delta t$, $T$, and $n$ parameters as defined in subsection III.A are 1 day, 10 days, and 2 respectively. Due to the discrete nature of the time-expanded network, there are only a few dates at which the optimizer can dispatch servicer $v_0$ to the satellite in need. In this project, it was decided to offset the service window by the $n$ spaceflight time steps so that the servicer would have time to fly over the network between $t = 20$ and $t = 22$ and reach the satellite for the first service opportunity defined at $\tau = 22$ on Fig. 3. In this work, the service opportunities (depicted by the red dots on Fig. 3) are defined at the last of the spaceflight time steps encompassed within the service window. Each of these service opportunities is then given a service assignment variable whose value is decided by the optimizer. The time set $\mathcal{W}_{s_0}$ associated with service need $s_0$ would then be $\{22,32\}$. This example illustrates faithfully what is happening under the hood of the proposed optimization framework and can easily be generalized to more complex situations.

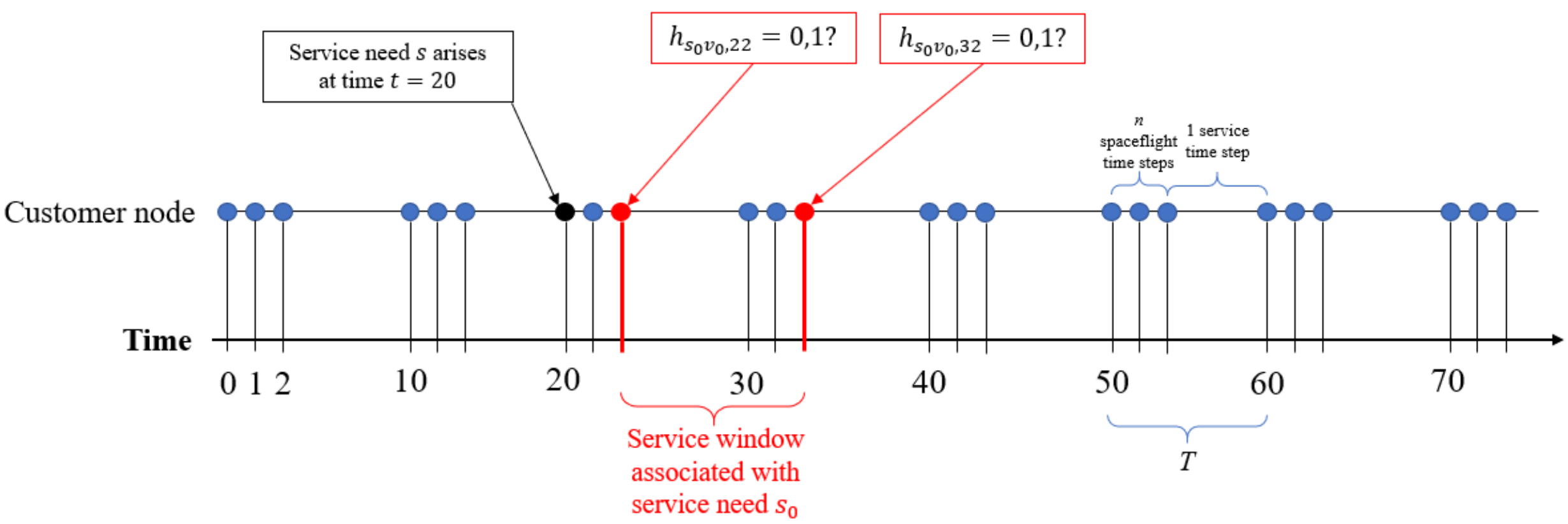

**Fig. 3 Illustration of the concept of service window and its relationship with the service assignment variables $h_{sv\tau}$ ($\Delta t$ = 1 day, $T$ = 10 days, and $n$ = 2).**

*Servicer dispatch variables*

The service assignment variables previously introduced cannot be used to enforce the presence of a servicer at a given node since they are defined for subsets only (*i.e.,* $\mathcal{W}_s$) of the entire set of time steps $\mathcal{T}$. To remedy this, the servicer dispatch variables, denoted $b_{svt}$, are defined for each time step in $\mathcal{T}$. Note that the time index of these variables is different from that of the service assignment variables ($t$ instead of $\tau$) to indicate that they are indexed over $\mathcal{T}$ and not over $\mathcal{W}_s$. At any time step $t$ where $b_{svt} = 1$, servicer $v$ must be in the process of providing service $s$. The optimizer sets it to 0 otherwise.

The servicer dispatch variables are related to the service assignment variables through the equality constraint:

$$b_{svt} = \sum_{\tau \in \mathcal{W}_s} h_{sv\tau} \beta_{s\tau t} \qquad \forall v \in \mathcal{V}_{\text{serv}}, \forall s \in \mathcal{S}, \forall t \in \mathcal{T} \tag{9}$$

The binary parameters $\beta_{s\tau t}$ defined in Eq. (9) are automatically generated by the OOS optimization framework before a planning horizon optimization based on the input data related to the services. Note that these parameters are indexed by both $t \in \mathcal{T}$ and $\tau \in \mathcal{W}_s$. Essentially, the binary parameter $\beta_{s\tau t}$ captures the duration of service $s$ and at which time steps $t$ the service has to be provided if the OOS operator decides to start the service at the date $\tau$.

The above example used to describe the service assignment variables (*cf* Fig. 3) is extended here to give more insight into Eq. (9). From the situation illustrated in Fig. 3, Eq. (9) would be specialized into Eq. (10). The notations from Fig. 3 are kept for consistency. Also note that $\mathcal{W}_{s_0}$ was replaced with {22,32} in Eq. (10).

$$b_{s_0 v_0 t} = \sum_{\tau \in \{22,32\}} h_{s_0 v_0 \tau} \beta_{s_0 \tau t} = h_{s_0 v_0,22} \beta_{s_0,22,t} + h_{s_0 v_0,32} \beta_{s_0,32,t} \qquad \forall t \in \mathcal{T} \tag{10}$$

Next, Table 8 illustrates what would be in this notional example the sequence of 1's and 0's for the $\beta_{s\tau t}$ parameters. In order to generate these parameters, it is assumed that the refueling service $s_0$ lasts for 40 days. Note that the 1's in Table 8 do not span 40 days; instead, this 40-day service is interpreted as over four service time steps (from after day 22 through right before day 60). This modeling is used to have the end of a service period match with the beginning of a service time interval $T$ so that, when the servicer is done providing a service, it can benefit from the transportation arcs defined over the $n$ spaceflight time steps (*cf* subsection III.A for the definition of the $T$ and $n$ parameters).

**Table 8. Notional sequence of binary parameters $\beta_{s\tau t}$ used to relate the $b_{svt}$ variable with the $h_{sv\tau}$ variables.**

| Time [day] | 0 | 1 | 2 | 10 | 11 | 12 | 20 | 21 | 22 | 30 | 31 | 32 | 40 | 41 | 42 | 50 | 51 | 52 | 60 | 61 | 62 | 70 | 71 | 72 |
|---|---|---|---|---|---|---|---|---|---|---|---|---|---|---|---|---|---|---|---|---|---|---|---|---|
| $\beta_{s_0,22,t}$ | 0 | 0 | 0 | 0 | 0 | 0 | 0 | 0 | **1** | **1** | **1** | **1** | **1** | **1** | **1** | **1** | **1** | **1** | 0 | 0 | 0 | 0 | 0 | 0 |
| $\beta_{s_0,32,t}$ | 0 | 0 | 0 | 0 | 0 | 0 | 0 | 0 | 0 | 0 | 0 | **1** | **1** | **1** | **1** | **1** | **1** | **1** | **1** | **1** | **1** | 0 | 0 | 0 |

Through the service assignment variables $h_{sv\tau}$, the optimizer then "turns on" the sequence of binary parameters that is optimal. The selected sequence of binary parameters is then assigned to the servicer dispatch variables $b_{svt}$. Figure 4 illustrates the selection of the sequence of binary parameters through $h_{sv\tau}$. Note that the optimizer is left with the choice to not provide a service which means that the $h_{sv\tau}$ variables in Eq. (9) may all be set to 0.

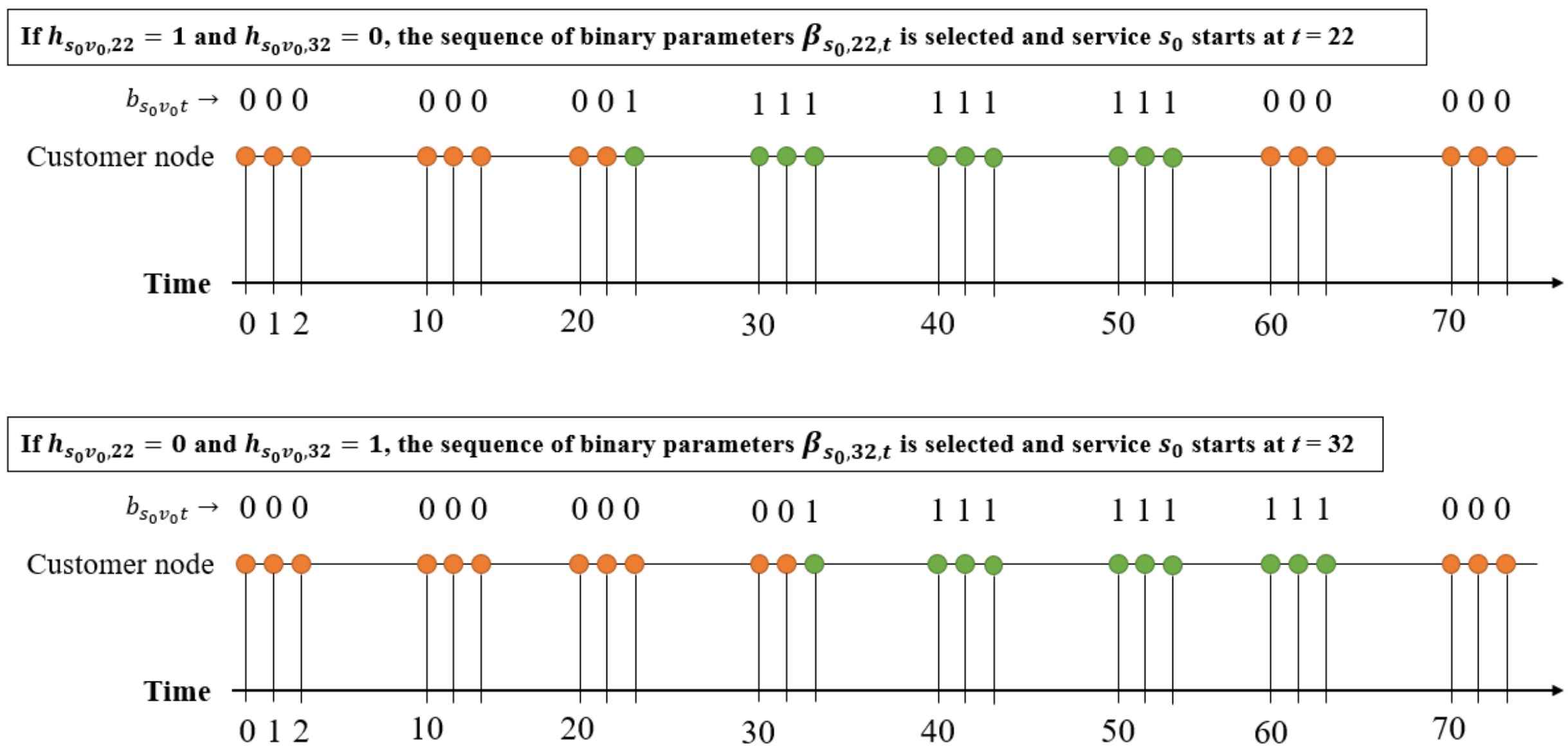

**Fig. 4 Illustration of how the servicer dispatch variables $b_{svt}$ are assigned a value by selecting a sequence of binary parameters $\beta_{s\tau t}$ through the service assignment variables $h_{sv\tau}$.**

### ii. OOS-specific constraints

The OOS-specific variables introduced in the previous section are used to define three additional constraints that are appended to the traditional space logistics formulation. The first constraint given by Eq. (11) specifies that each service $s$ may be scheduled at most once.

$$\sum_{v \in \mathcal{V}_{\text{serv}}} \sum_{\tau \in \mathcal{W}_s} h_{sv\tau} \leq 1 \qquad \forall s \in \mathcal{S} \qquad (11)$$

The second constraint given by Eq. (12) specifies that at most one service may be provided to a customer satellite per time step. This means that if a customer satellite experiences a refueling service need and a repair service need at the same time, the OOS infrastructure will provide at most one of them. This constraint essentially prevents the clustering of servicers at a customer satellite in order to ensure the safest possible operations. The index set $\mathcal{N}_c$ in Eq. (12) is a subset of the set of nodes $\mathcal{N}$ that includes only the customer nodes.

$$\sum_{v \in \mathcal{V}_{\text{serv}}} \left\{ \sum_{s \in \mathcal{S}_i} b_{svt} \right\} \leq 1 \qquad \forall i \in \mathcal{N}_c , \forall t \in \mathcal{T} \tag{12}$$

The final constraint given by Eq. (13) is critical for an accurate modeling of the OOS operations as it forces servicers at customer nodes for the entire durations of the services. This is done through the definition of a set of binary parameters $\gamma_{sk}$ that map a tool $k$ to service $s$ as given in Table 5. For example, if service $s$ is a refueling service, then $\gamma_{sk}$ = 1 for $k$ being the refueling apparatus tool whereas $\gamma_{sk}$ = 0 otherwise. The index set of the tool commodities is denoted $\mathcal{K}_{\text{tools}}$ and is a subset of the set $\mathcal{K}_I$ of non-vehicle integer commodities defined in section III.C.1. The servicers are forced at customer nodes by requiring the use of their tools for services.

$$x^{+}_{viitk} \geq \sum_{s \in \mathcal{S}_i} \gamma_{sk} b_{svt} \qquad \forall v \in \mathcal{V}_{\text{serv}}, \forall i \in \mathcal{N}_c, \forall t \in \mathcal{T}, \forall k \in \mathcal{K}_{\text{tools}} \tag{13}$$

**iii. OOS-specific objective function**

Although most literature in space logistics has focused on the launch mass as the objective function, different objective functions may be defined based on the true motivations of the OOS venture. This paper assumes that the modeled OOS company wishes to maximize the profits over a planning horizon. The profits are defined as the difference between the revenues generated by the provision of the services and the costs incurred by the infrastructure.

The revenues are captured through the service assignment variables $h_{sv\tau}$ and revenue parameters $r_s$ defined for each service. The revenues are modeled by Eq. (14).

$$J_{\text{revenues}} = \sum_{v\in\mathcal{V}_{\text{serv}}} \sum_{i\in\mathcal{N}_c} \sum_{s\in\mathcal{S}_i} \sum_{\tau\in\mathcal{W}_s} \{r_s h_{sv\tau}\} \tag{14}$$

The costs include five terms: (a) the purchase, development and manufacturing costs $J_{\text{pdm}}$ of the commodities, modeled by Eq. (15); (b) the launch costs $J_{\text{launch}}$ of the commodities sent into space from Earth, modeled by Eq. (16); (c) the penalty fees for delayed services $J_{\text{delay}}$, modeled by Eq. (17); (d) the operating costs of the depots $J_{\text{depots}}$, modeled by Eq. (18); and (e) the operating costs of the servicers $J_{\text{serv}}$, modeled by Eq. (19). The meaning of the parameters and index sets involved in the equations can be found in the nomenclature.

$$J_{\text{pdm}} = \sum_{v\in\mathcal{V}} \sum_{t\in\mathcal{T}} \sum_{i\in\mathcal{N}_e} \sum_{j\in\mathcal{N}} \left\{ \sum_{k\in\mathcal{K}_C} c_k^{\text{pdm}} x_{vijtk}^{+} + \sum_{k\in\mathcal{K}_I} c_k^{\text{pdm}} m_k x_{vijtk}^{+} + c_v^{\text{pdm}} y_{vijt}^{+} \right\} \tag{15}$$

$$J_{\text{launch}} = c^{\text{launch}} \sum_{v\in\mathcal{V}} \sum_{t\in\mathcal{T}} \sum_{i\in\mathcal{N}_e} \sum_{j\in\mathcal{N}} \left\{ \sum_{k\in\mathcal{K}_C} x_{vijtk}^{+} + \sum_{k\in\mathcal{K}_I} m_k x_{vijtk}^{+} + m_v y_{vijt}^{+} \right\} \tag{16}$$

$$J_{\text{delay}} = \sum_{i\in\mathcal{N}_c} \sum_{s\in\mathcal{S}_i} \sum_{\tau\in\mathcal{W}_s} \left\{ c_s^{\text{delay}} (\tau - \tau_s) \sum_{v\in\mathcal{V}_{\text{serv}}} h_{sv\tau} \right\} \tag{17}$$

$$J_{\text{depots}} = \sum_{v\in\mathcal{V}_{\text{depot}}} \sum_{i\in\mathcal{N}_p} \sum_{t\in\mathcal{T}} \{ c_v^{\text{depot}} \Delta_{iit} y_{viit}^{+} \} \tag{18}$$

$$J_{\text{serv}} = \sum_{v\in\mathcal{V}_{\text{serv}}} \sum_{i,j\in\mathcal{N}_c\cup\mathcal{N}_p} \sum_{t\in\mathcal{T}} \{ c_v^{\text{serv}} \Delta_{ijt} y_{vijt}^{+} \} \tag{19}$$

The profits to be maximized are then simply modeled by Eq. (20).

$$J = J_{\text{revenues}} - J_{\text{pdm}} - J_{\text{launch}} - J_{\text{delay}} - J_{\text{depots}} - J_{\text{serv}} \tag{20}$$

### C. Application of Rolling Horizon Decision Making to On-Orbit Servicing

In order for an OOS venture to make long-term and sustainable decisions early on, it must account for uncertainties in the servicing demand. In this paper in particular, uncertainty stems from the randomly-occurring service needs. Efficient methods such as the RH approach can be leveraged to support the decision making of OOS operators.

RH decision making is used to make decisions in a dynamic stochastic environment, typically characterized by uncertainties in demand and the resulting cost in forecasting this demand. More specifically, this technique consists in making optimal a series of decisions repeatedly based on short- to mid-term forecast of future demand. An advantage of RH is to decompose a complex dynamic scheduling problem into smaller sub-problems whose combined optimal solutions yield a satisfying solution of the original problem at a lesser computational cost. Readers are referred to Refs. [25,26] for a more complete theoretical background on RH decision making.

The RH approach consists in dividing the total time horizon (called *Scheduling Horizon* (SH) in this paper) in shorter time intervals (called *Planning Horizons* (PH)) over which the operations of a given system are optimized assuming an accurate demand forecast in those intervals. The results of the optimization over the first period of each PH (called *Control Horizon* (CH) in this paper) are retained before "rolling on" the PH over time and reiterating on the described procedure. A critical step in the RH algorithm is to ensure the continuity of the state of the system between the end of a CH and the beginning of the following PH. Figure 5 concisely summarizes the traditional RH procedure with three successive PHs and CHs.

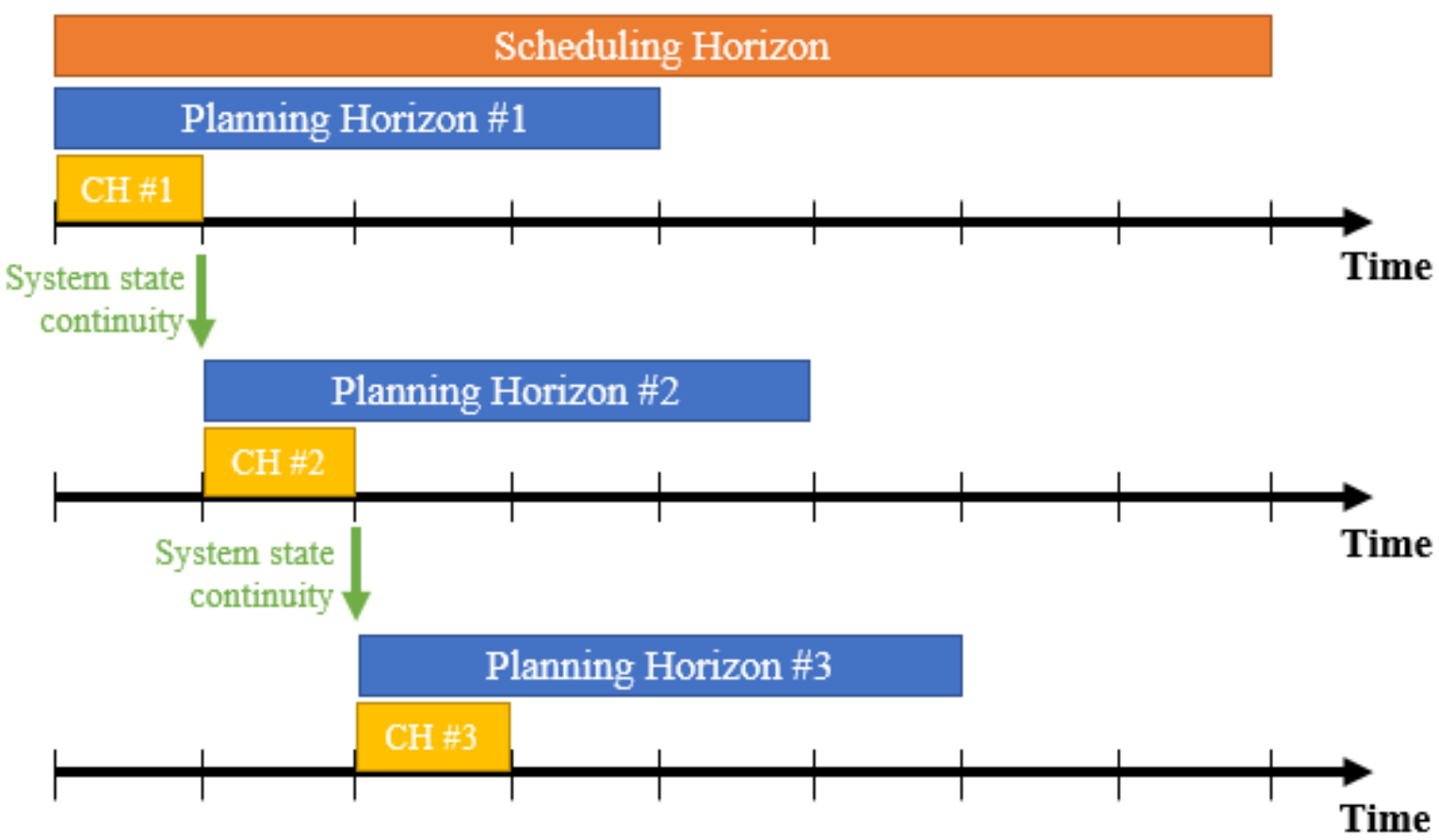

**Fig. 5 An illustration of the traditional Rolling Horizon procedure.**

As depicted in Fig. 5, the traditional RH approach re-plans the operations of the system every fixed time interval. This is possible when the uncertain demand can be forecast over the PH. In the context of OOS however, predicting the occurrences of random service needs such as satellite repairs is a difficult task and not appropriate for realistic operational scheduling. Therefore, a modification of the traditional RH procedure is proposed in which the OOS operations are re-planned whenever a random service need occurs. The OOS operations may also be re-planned on a regular basis if there has not been any random service need for a long time. The demand forecast over each PH thus comprises the random service need that triggers the re-planning of the operations (if any) and the deterministic service needs arising over the PH. Finally, as for the traditional RH procedure, the state of the OOS infrastructure must be propagated from one CH to the next PH. Here, the state at any date $t$ of the OOS infrastructure comprises the position and amount of all commodities in the static network at $t$.

Figures 6, 7, 8 illustrate the proposed RH procedure in the context of OOS operational scheduling for a notional 6-node dynamic network. On those figures, the detailed time-expansion of the static network as presented in Fig. 2 is omitted for readability. A notional service demand history is given to keep track over time of the service needs, their occurrence dates, and to which customer satellite the services must be provided. Figure 6 illustrates a notional schedule of a servicer over the first PH. The servicer's path over the dynamic network is represented by bold, yellow direct arcs. The operations are automatically re-planned after 120 days of operations as no random service need arises over that time interval (the 120-day time interval for automatic re-planning is user-defined). Figure 7 illustrates the servicer's schedule after the PH is rolled over between the dates $t = 120$ and $t = 480$. Note how the servicer's schedule from $t = 0$ to $t = 120$ is retained as part of the final solution whereas the servicer's schedule is modified after the date $t = 120$ due to a new deterministic service need, called $d_5$, in the forecast. At the date $t = 200$, a random service need called $r_1$ triggers a re-planning of the OOS operations. Thus, one can see that the CH in the proposed methodology does not have a fixed length contrary to the traditional RH procedure as shown in Fig. 5. Finally, Fig. 8 shows that the servicer's schedule is modified after the date $t = 200$ in order to provide a service in response to the random need $r_1$. Figures 6, 7, 8 show the path of only one servicer for readability, but the proposed methodology is flexible enough to consider as many servicers as desired.

The OOS operations over each PH of the RH procedure are optimized by formulating the OOS logistics problem as a MILP model presented in subsection III.B. Thus, the proposed methodology consists in solving a series of MILP problems over a given SH.

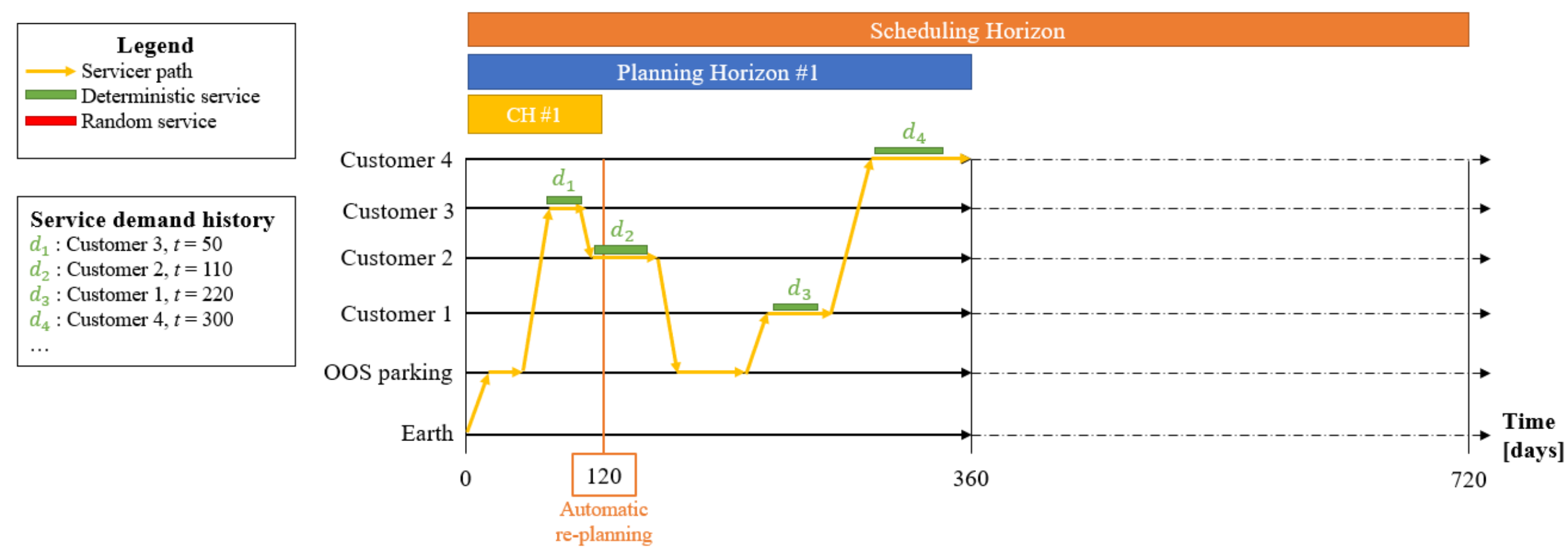

**Fig. 6 On-Orbit Servicing operations over the first planning and control horizons.**

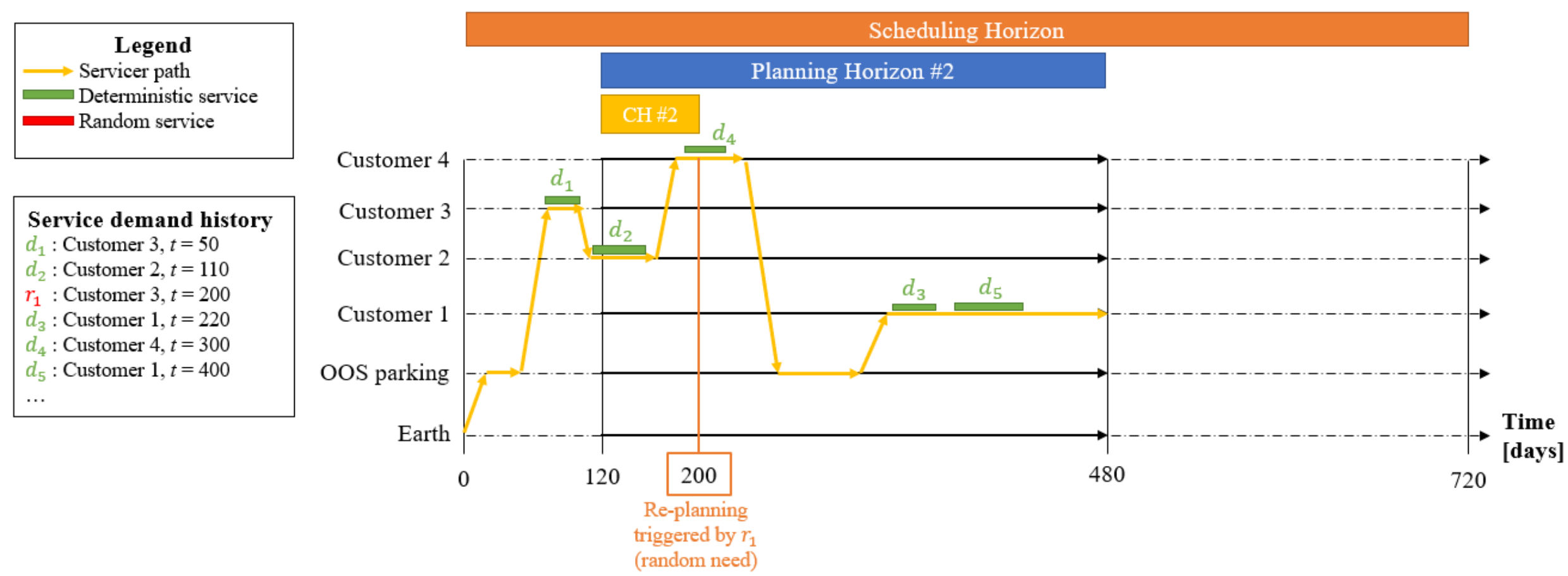

**Fig. 7 On-Orbit Servicing operations over the second planning and control horizons.**

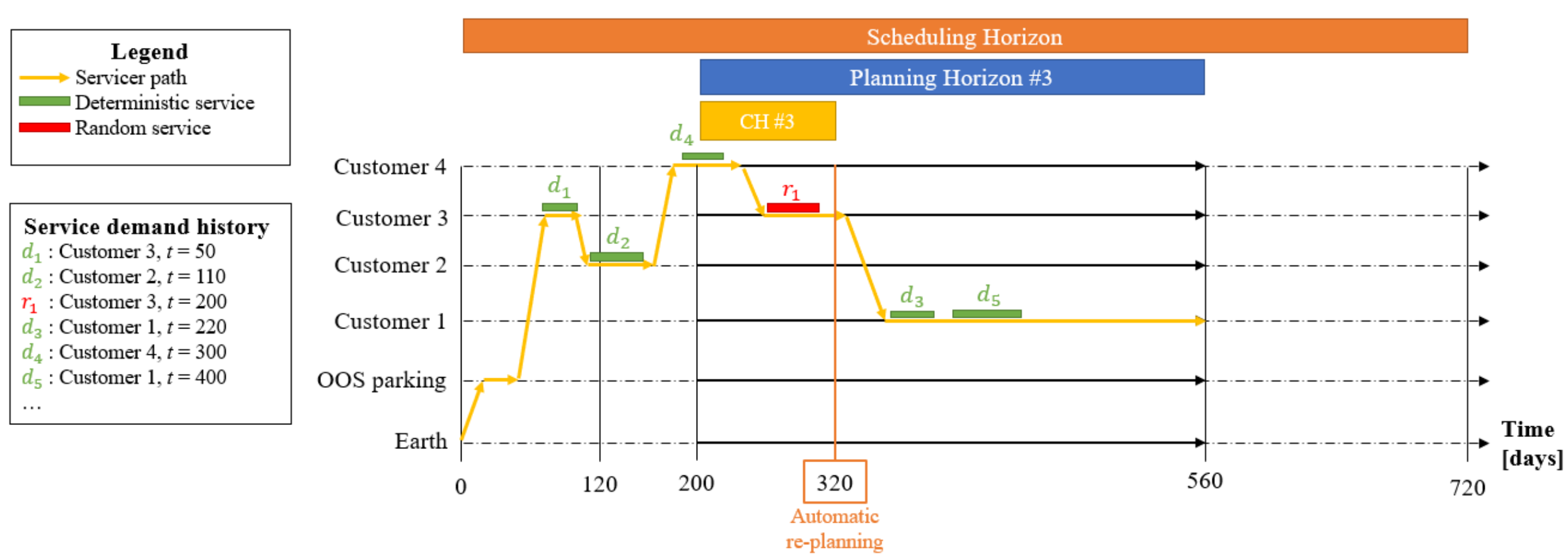

**Fig. 8 On-orbit Servicing operations over the third planning and control horizons.**

## IV. Case Studies

This section demonstrates two use cases of the OOS logistics optimization framework: a short-term operational scheduling use case; and a long-term strategic planning use case. Subsection IV.A summarizes the assumptions used to run the different scenarios designed to demonstrate the use cases. Subsection IV.B demonstrates the proposed method as a tool for short-term operational scheduling of an existing OOS infrastructure. Subsection IV.C demonstrates the proposed method as a tool to inform long-term strategic design of OOS architectures by comparing the performance of different OOS architectures.

### A. Assumptions and Scenarios

In this subsection, we first present the assumptions and data associated with the fleet of customer satellites. Then, the assumptions and data associated with the OOS infrastructure are given. Finally, the scenarios associated with the two use-cases are presented.

### i. Customer Fleet Assumptions

The data related to the deterministic and random service needs are given in Table 9 and Table 10 respectively. Note that the given data are for one customer satellite; by increasing the number of customer satellites in the simulations, the service need rates of the entire fleet of customer satellites increase. The positions of the customer satellites were retrieved from the UCS (Union of Concerned Scientists) database of satellites [28].

**Table 9. Assumptions related to the deterministic service needs (*).**

| | Inspection | Refueling | Station Keeping |
|---|---|---|---|
| **Revenues [$M]** | 10 [7] | 15 [7] | 20 |
| **Delay penalty fee [$/day]** | 5,000 [7] | 100,000 [7] | 100,000 [7] |
| **Service duration [days]** | 10 [7] | 30 [7] | 180 [14] |
| **Service window [days]** | 30 | 30 | 30 |
| **Frequency of occurrence [days]** | 6,310 [11] | 2,100 [11] | 2,100 [11] |

* References are indicated within brackets. Data without reference is assumed. The frequency of occurrence is derived from the data given in Ref. [11].

**Table 10. Assumptions related to the random service needs (*).**

| | Repositioning | Retirement | Repair | Mechanism Deployment |
|---|---|---|---|---|
| **Revenues [$M]** | 10 [7] | 10 [7] | 30 | 25 [7] |
| **Delay penalty fee [$/day]** | 100,000 [7] | 0 [7] | 100,000 [7] | 100,000 [7] |
| **Service duration [days]** | 30 [7] | 30 [7] | 30 [7] | 30 [7] |
| **Service window [days]** | 30 | 30 | 30 | 30 |
| **Mean frequency of occurrence [days]** | 2,520 [11] | 2,520 [11] | 9,020 [11] | 21,050 [11] |

* References are indicated within brackets. Data without reference is assumed. The mean frequency of occurrence is derived from the data given in Ref. [11].

### ii. OOS Infrastructure Assumptions

The four notional servicer tools given in Table 5 each have an assumed cost of $100,000 and an assumed mass of 100kg.

Two types of servicers are defined: 1 versatile servicer with 4 integrated tools and 4 specialized servicers each with exactly one tool of each type. The versatile servicer can provide all seven defined services whereas the specialized servicers can only provide the services their tools are suited for. The detailed assumptions are given in Table 11. The baseline dry mass for the versatile servicer is taken from Ref. [7] but then reduced for the specialized servicers which are assumed to be less capable so of smaller size. The baseline propellant capacity in Ref. [7] is 500kg, but is increased in this paper because we are modeling sustainable infrastructures in this project. Thus, the servicers must be able to return to the depot to get refueled. The assumed specific impulse is characteristic of bipropellant rocket engine. Furthermore, this paper assumes the refueling of the servicers to be instantaneous operations. The justification behind this assumption is that, as OOS operations become routine, refueling of the servicers will likely not take more than one time step (*i.e.*, 2 days) [30]. This assumption can be modified depending on the technology performance.

**Table 11. Assumptions related to the servicers (*).**

| | **Versatile (V)** | **Specialized 1 (S1)** | **Specialized 2 (S2)** | **Specialized 3 (S3)** | **Specialized 4 (S4)** |
|---|---|---|---|---|---|
| **Tools** | T1, T2, T3, T4 | T1 | T2 | T3 | T4 |
| **Dry mass [kg]** | 3,000 [7] | 2,000 | 2,000 | 2,000 | 2,000 |
| **Propellant capacity [kg]** | 1,000 | 1,000 | 1,000 | 1,000 | 1,000 |
| **Specific Impulse [s]** | 316 | 316 | 316 | 316 | 316 |
| **Manufacturing cost [$M]** | 75 | 50 [29] | 50 [29] | 50 [29] | 50 [29] |
| **Operating cost [$/day]** | 13,000 [15] | 13,000 [15] | 13,000 [15] | 13,000 [15] | 13,000 [15] |

* References are indicated within brackets. Data without reference is assumed.

An orbital depot is assumed pre-deployed at a GEO orbital slot located at a longitude of 170 deg West (over the Pacific Ocean). The depot is assumed to consume its own monopropellant at a rate of 0.14kg/day for station keeping [31]. The manufacturing and operating costs of the depot are assumed to be $200M and $13,000/day respectively.

The launch vehicle used for this analysis is based off of a Falcon 9 launcher. One launcher is assumed to launch every 30 days to resupply the depot with a maximum payload capacity of 8,300kg. The launch cost is assumed to be $11,300/kg.

Finally, the commodities considered in the case studies are the spares (assumed price tag: $1,000/kg), bipropellant for the servicers (price tag for monomethyl hydrazine: $180/kg), and monopropellant (price tag for hydrazine: $230/kg).

### iii. Use-Cases' Scenarios

Using the assumptions presented previously, several scenarios are designed to demonstrate the short-term operational scheduling and long-term strategic planning capabilities of the proposed methodology. Three different levels of service demand are defined by considering 30 satellites (low demand rate), 71 satellites (medium demand rate), or 142 satellites (high demand rate).

The first case study consists in demonstrating the value of the proposed framework in optimizing short-term operational scheduling. This situation occurs when an OOS operator wishes to optimally re-plan the operations of its infrastructure for the next few months given a forecast of service needs. The OOS infrastructure is assumed to be initially pre-deployed at some customer satellites to account for the fact that at re-planning some servicers may already be providing services. Note that the proposed framework allows users to define any initial state of the OOS infrastructure. For the purpose of this demonstration, the framework is run for a single PH optimization. This scenario is summarized in Table 12.

**Table 12. Scenario definition for the use case of short-term operational scheduling.**

| | |
|---|---|
| **Servicers** | 4 specialized servicers |
| **Depot** | 1 depot pre-deployed at 170deg West |
| **Planning horizon** | 100 days |
| **Customer fleet** | 142 customer satellites (high demand) |
| **Initial conditions** | - S1 pre-deployed at satellite *UFO-4* providing refueling service<br>- S2 pre-deployed at satellite *UFO-4* (due to prior service)<br>- S3 pre-deployed at satellite *Echostar 23* (due to prior service)<br>- S4 pre-deployed at satellite *Galaxy 28* providing station keeping service |

The second case study demonstrates the value of the proposed framework in analyzing the pros and cons of a particular OOS architecture over the long term and trading different OOS architectures with one another. In order to assess the value of each OOS architecture, the optimization framework is run for a scheduling horizon of 5 years. Each planning horizon optimization is performed over 100 days and the length of the control horizon in case of automatic re-planning is 60 days. For this case study, we consider two different OOS architectures: a monolithic architecture using 1 versatile servicer; and a distributed architecture using 4 specialized servicers. The value and profitability of each architecture is assessed over 5 years of operations under low service demand (30 satellites), medium service demand (71 satellites) and high service demand (142 satellites). The value of an OOS architecture is defined as the difference between the profits and the initial investments. Thus, an OOS architecture is considered valuable at a given date if its value is positive. Table 13 summarizes the parameters related to each OOS architecture.

**Table 13. Definition of the OOS architectures for use case of long-term strategic planning.**

| | **Monolithic architecture** | **Distributed architecture** |
|---|---|---|
| **Servicers** | 1 versatile | 4 specialized |
| **Depot** | 1 depot pre-deployed at 170deg West | 1 depot pre-deployed at 170deg West |
| **Planning horizon** | 100 days | 100 days |
| **Control horizon (in case of automatic re-planning)** | 60 days | 60 days |
| **Initial conditions** | The servicer is pre-deployed at the OOS parking node | All servicers are pre-deployed at the OOS parking node |

### B. Use Case 1: Short-Term Operational Scheduling

Assume that an OOS company owns a fleet of vehicles to provide orbital services to a set of satellites all located in GEO. Whenever it decides to re-plan its operations, the OOS optimization framework would be run for a few hours (depending on the problem size) and gets an output of the flow of operations. Note that not all customer satellites will need a service during the period considered for re-planning. The framework thus includes in the network only the satellites that experience service needs over the PH. This minimizes the size of the optimization problem thus allowing for a faster solving. For example, in the scenario presently considered, the customer base is of 142 satellites. However, over the considered PH, only 10 satellites are included in the network as the other satellites do not experience any service need during that period. Of these 10 satellites, only 8 are visited in the final solution. Figure 9 gives the network for the considered scenario. The customer satellites are represented by colored dots. The red dots correspond to the visited customer satellites, whereas the blue dots are not visited by any servicer.

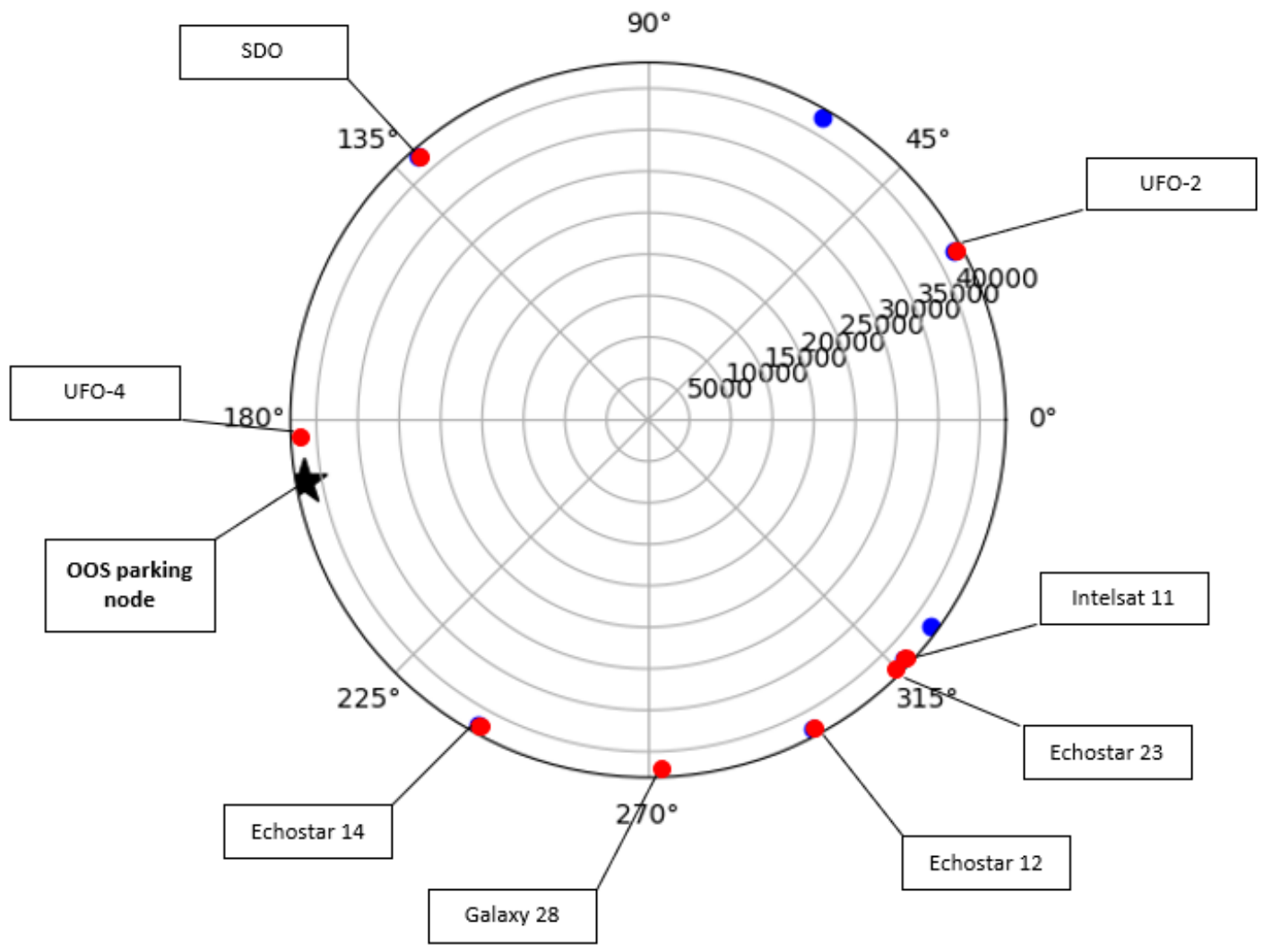

**Fig. 9 Static network used to demonstrate the short-term operational scheduling of the OOS infrastructure.**

The optimization framework is run on an Intel® Core™ i7-9700, 3.00GHz platform with the Gurobi 9 optimizer. The solution was reached in less than 10 seconds with a gap of 1% as stopping criterion. For this case study, the time-expansion parameters are $\Delta t = 2$ days, $T = 10$ days, $n = 2$. The outputs of the optimization are depicted on Figs. 10 and 11.

It is worth noting on Figs. 10 and 11 that the framework optimally leverages the service windows defined to give the OOS planners more flexibility in assigning the services to the servicers. For example, the S2 servicer is assigned the "Inspection" service at the *Echostar 12* satellite before that at the *UFO-2* satellite whereas the service need at *UFO-2* arises first. The optimizer decides to plan this way the operations of servicer S2 because *Echostar 12* is "on the way" between the depot and *UFO-2*. Doing so, servicer S4 saves propellant and thus limits the cost associated with its operations. The same observation can be made between the "Refueling' services provided at *Echostar 14* and *SDO*.

Also noticeable is the return of servicer S1 to the depot at times 22 and 92. The servicer returns to the depot and then immediately departs for other services. This is because the servicer needs to be refueled by the depot, which is assumed to be instantaneous as mentioned earlier.

For information, the optimal operations presented on Figs. 10 and 11 generate based on our initial assumptions \$65M of revenues for a total operating cost of \$46.5M. As depicted on the figures, the operations of the OOS infrastructure are profitable. These values however may not be final for the period $t = 0\ldots100$ as re-planning may be needed at a later time due to a randomly occurring service need. When a random need occurs, the OOS planners would re-run the framework with the initial conditions corresponding to the current orbital state of the OOS infrastructure and an updated set of service needs.

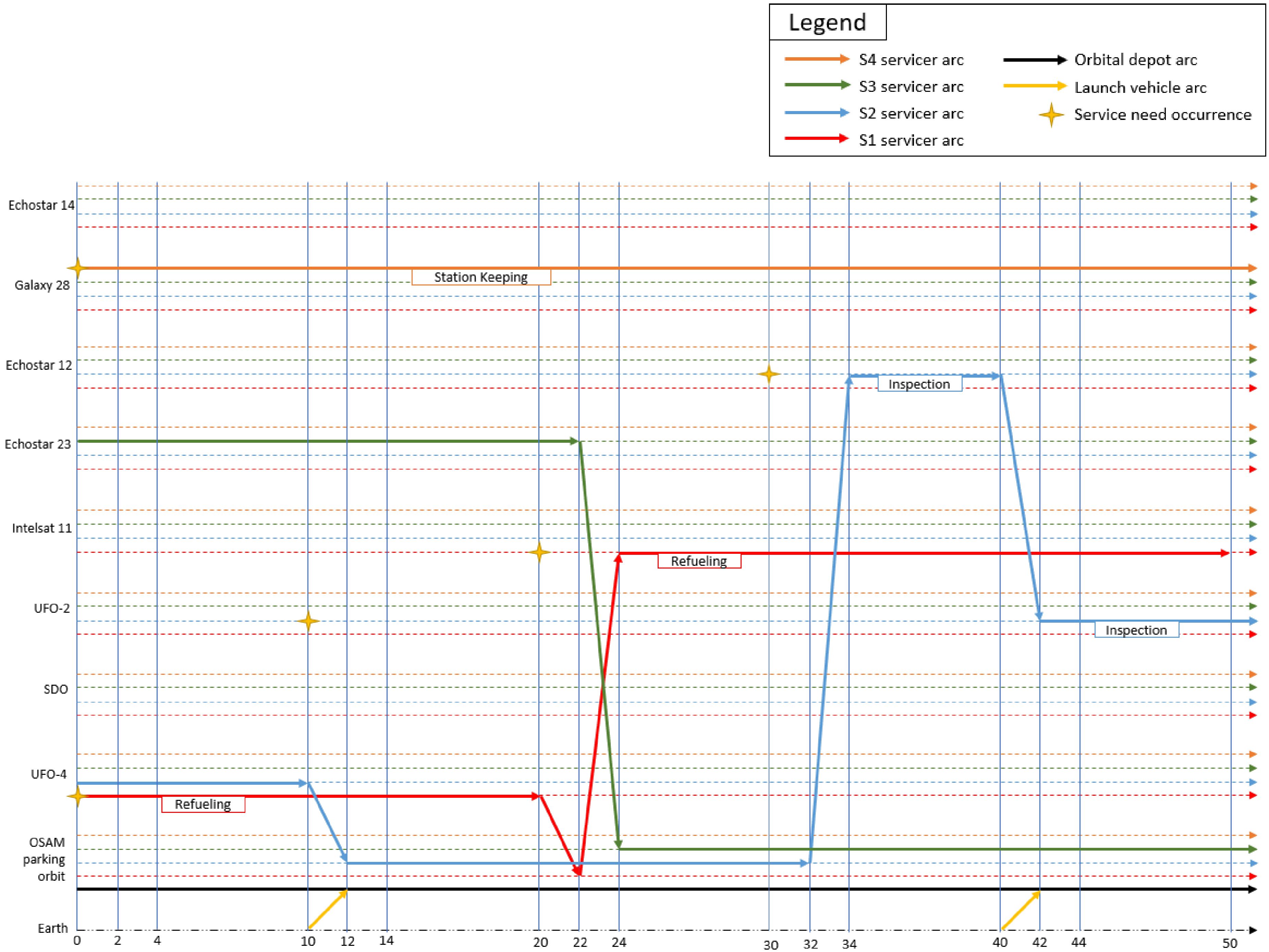

**Figure 10. Breakdown over the first 50 days of the operational scheduling of a notional OOS infrastructure consisting of 4 specialized servicers and a depot ($\Delta t$ = 2 days, $T$ = 10 days, $n$ = 2).**

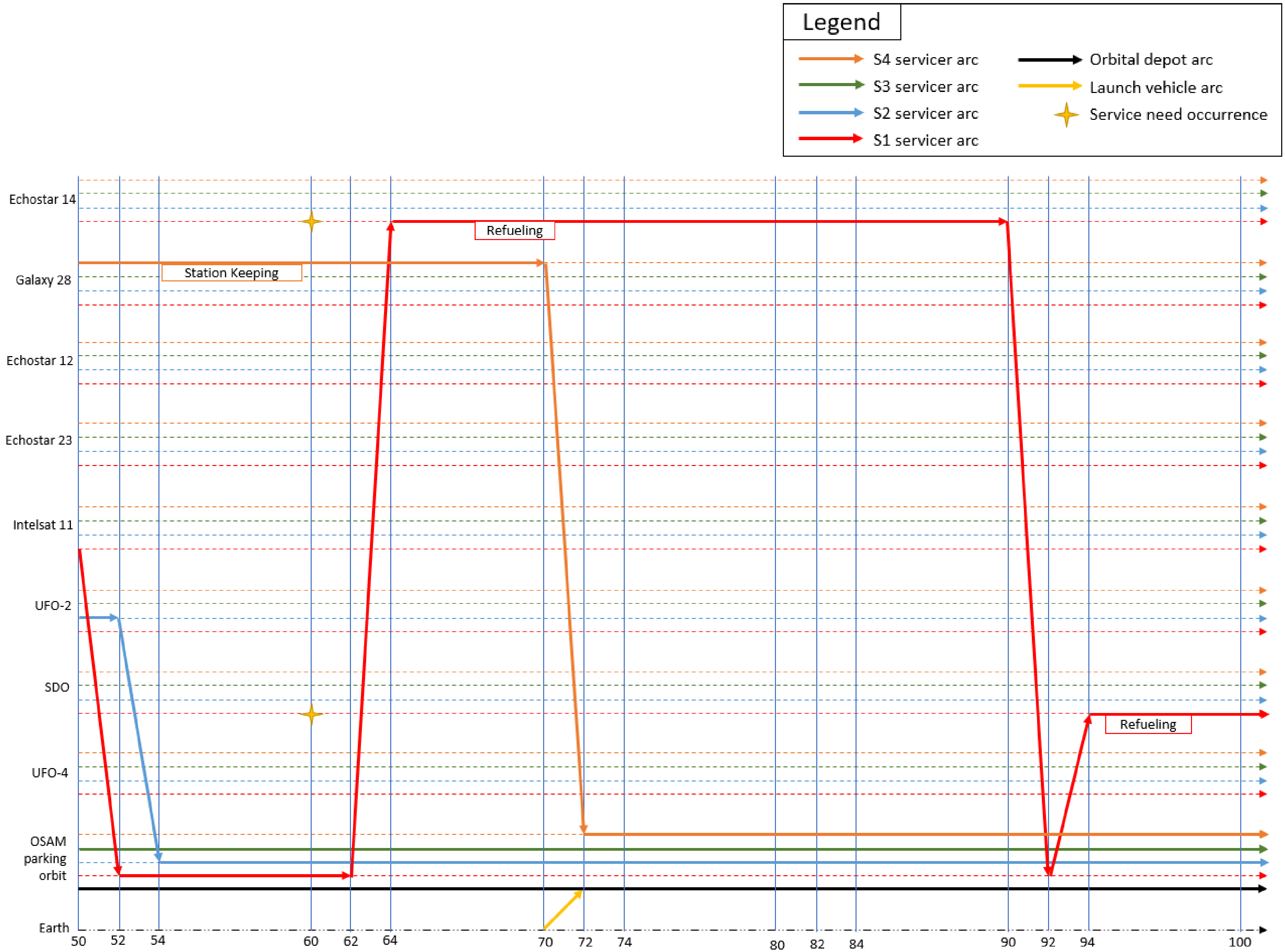

**Figure 11. Breakdown over the last 50 days of the operational scheduling of a notional OOS infrastructure consisting of 4 specialized servicers and a depot ($\Delta t$ = 2 days, $T$ = 10 days, $n$ = 2).**

## C. Use Case 2: Long-Term Strategic Planning

Besides being used as a short-term operational scheduling tool, the optimization framework may also be leveraged as a long-term strategic planning tool to assess the advantage of a given OOS architecture over another one and/or assess the impact of exogenous variables (*i.e.*, non-controllable by the OOS planners) such as the state of the on-orbit servicing market (*e.g.,* service demand level). In this paper, the OOS architectures are compared based on the profits they generate and their values at the end of the scheduling horizon. This case therefore demonstrates the developed method as a tool to design the architecture for OOS considering a long time horizon.

There are some additional considerations that need to be considered for the long-term strategic planning. To plan the operations of an OOS infrastructure over the short-term, only deterministic service needs are considered. These include for instance the services resulting from a contract between the OOS operator and a GEO satellite operator who wants its satellites to be regularly refueled. In the long-term strategic planning context, the uncertainty associated with services that cannot be planned in advance (*e.g.,* a repair) is captured through the RH approach. Two types of service needs are random by nature: repair and mechanism deployment. To some extent, retirement and repositioning service needs may also be considered random, as assumed here.

As a demonstration example, the value and profitability of a distributed architecture (*i.e.,* 4 small servicers carrying a different tool) are compared to that of a monolithic architecture (*i.e.,* 1 large servicer carrying all tools needed to provide any service). OOS planners may also want to see the impact of the market on the metrics for each architecture. For example, a costly architecture such as one involving many small, specialized servicers may not be as profitable as a monolithic one if the demand is low. But what if the market is expected to significantly grow in the near term? Wouldn't a distributed architecture be more fit to capture the extra sources of revenue and thus outweigh its initial investment at a faster rate? Those are questions that the framework can help answer.

In order to demonstrate the value of the framework for long-term strategic planning, it is run for the two architectures outlined in Table 13 over a scheduling horizon of 5 years and for three levels of service demand (30, 71 and 142 satellites). The successive optimizations were solved on an Intel® Core™ i7-9700, 3.00GHz platform with the Gurobi 9 optimizer. The two stopping criteria used in the optimizations were a gap of 1% and a time limit of 7200s. In total, 104 successive optimizations were solved in an average time of 690s and with an average gap of 1.22%. Four optimizations reached the 7200s time limit before reaching the required gap; in which case the best feasible solutions are used. As for the previous case study, the time-expansion parameters are $\Delta t$ = 2 days, T = 10 days, and $n$ = 2.

We first assume that the customer base of the OOS infrastructure is of 30 satellites. Figure 12 compares the evolution of the values of the monolithic and distributed architectures. Figures 13 and 14 respectively represent the revenues generated and costs incurred by both infrastructures over time. The architectures are equivalent in terms of profitability as shown by the slopes of the value curves. In this case, the single servicer in the monolithic architecture is not overwhelmed by too large a number of service needs and can thus provide almost as many services as the distributed architecture with an overall lower cost. Note that based on the assumptions summarized in Table 11, the initial investment needed for the distributed architecture is significantly larger. This would need to be accounted for by decision makers regarding the design of an OOS infrastructure.

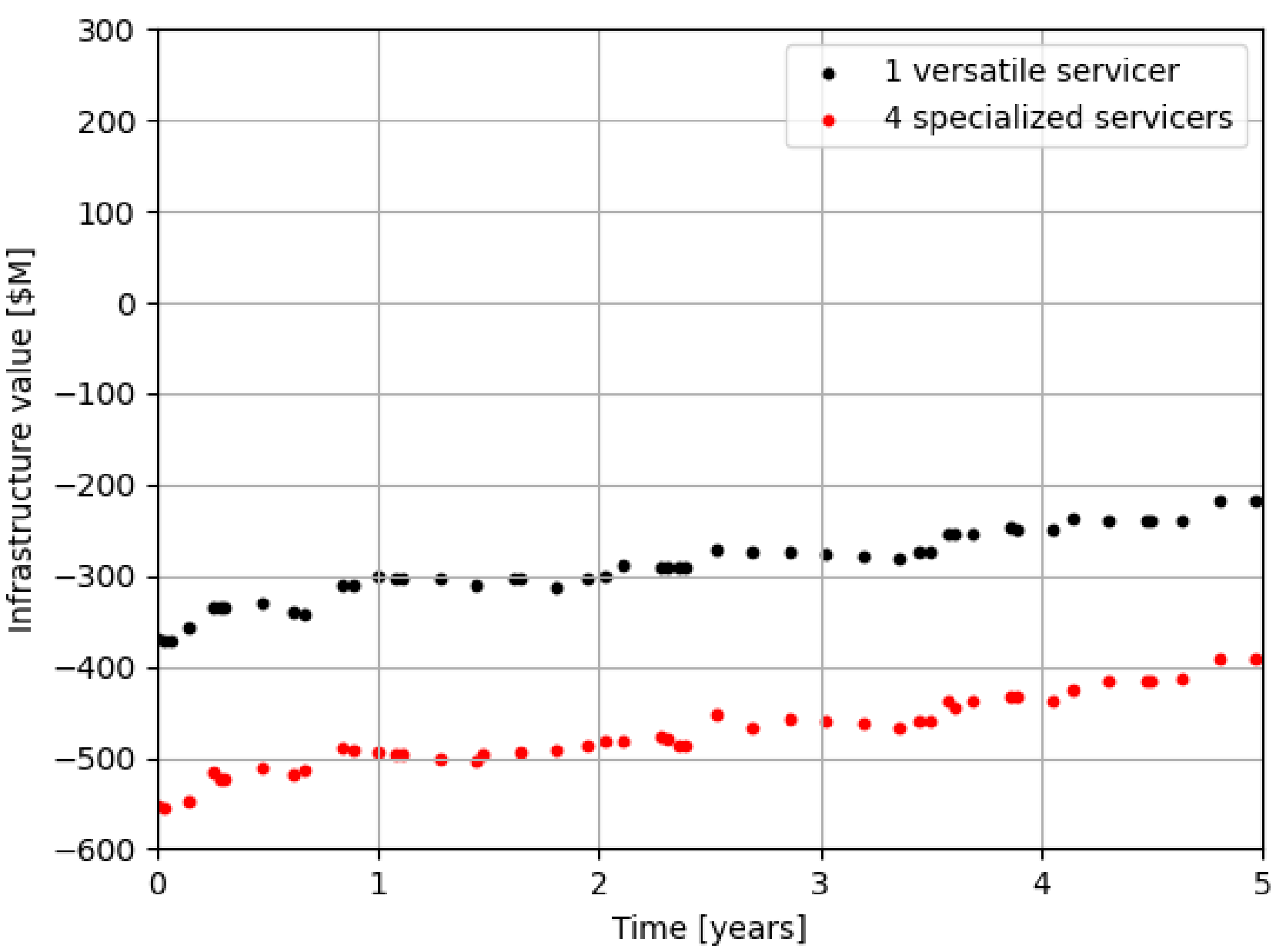

**Fig. 12 Values of the monolithic and distributed architectures for a customer base of 30 satellites.**

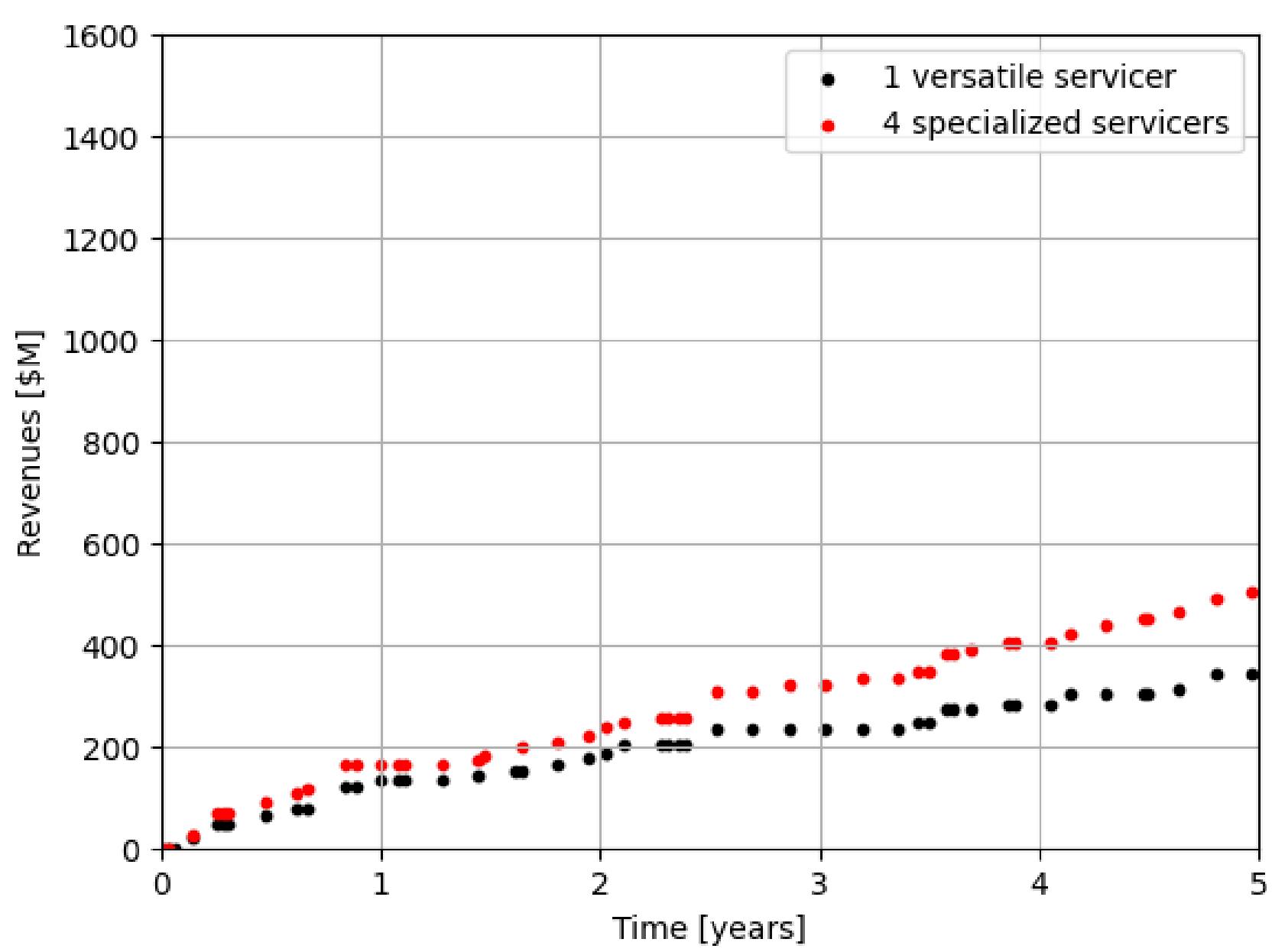

**Fig. 13 Revenues generated by the monolithic and distributed architectures for a customer base of 30 satellites.**

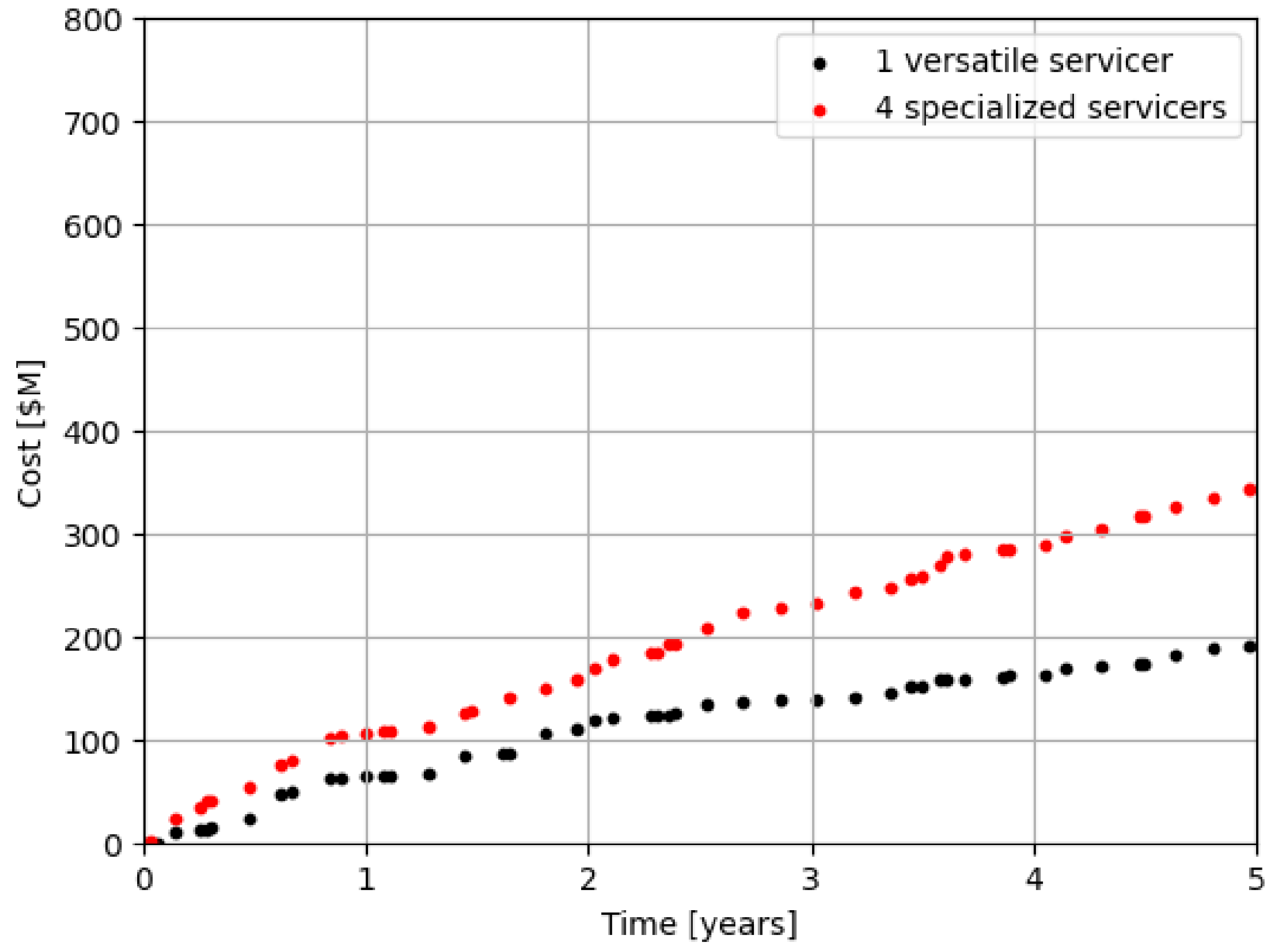

**Fig. 14 Costs incurred by the monolithic and distributed architectures for a customer base of 30 satellites.**

We next assume that the customer base more than doubles (71 satellites). Figure 15 represents the values of both architectures. Figures 16 and 17 respectively represent the revenues generated and costs incurred by both infrastructures over time. It is now clear that the distributed architecture outpaces the monolithic architecture in terms of profitability. It even becomes more valuable after 3 years of operation. In order to put Figs. 12 and 15 in perspective, assume that the forecast for the first 5 years is on the order of 30 satellites and that the following 5 years, the customer base doubles. Two of the decisions OOS planners may face are:

- wait until the market becomes more favorable so as to avoid a long payback period; or
- invest early in a distributed architecture whose profitability after a few years of operations would increase due to a growing customer base.

Investing early, the OOS operator may gain a leading competitive edge against new entrants as the market grows. This would also mean convincing investors to endure a longer payback period.

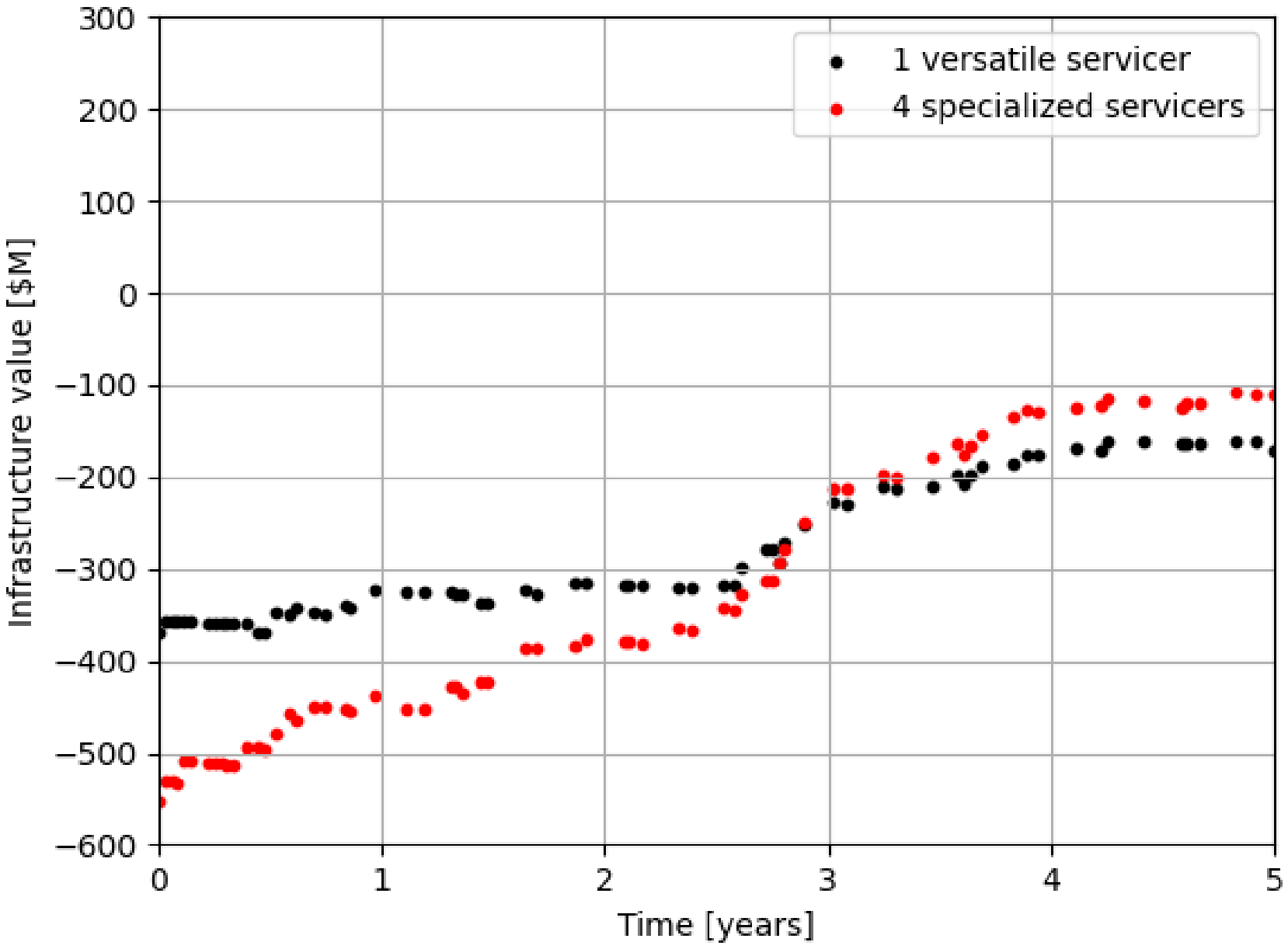

**Fig. 15 Values of the monolithic and distributed architectures for a customer base of 71 satellites.**

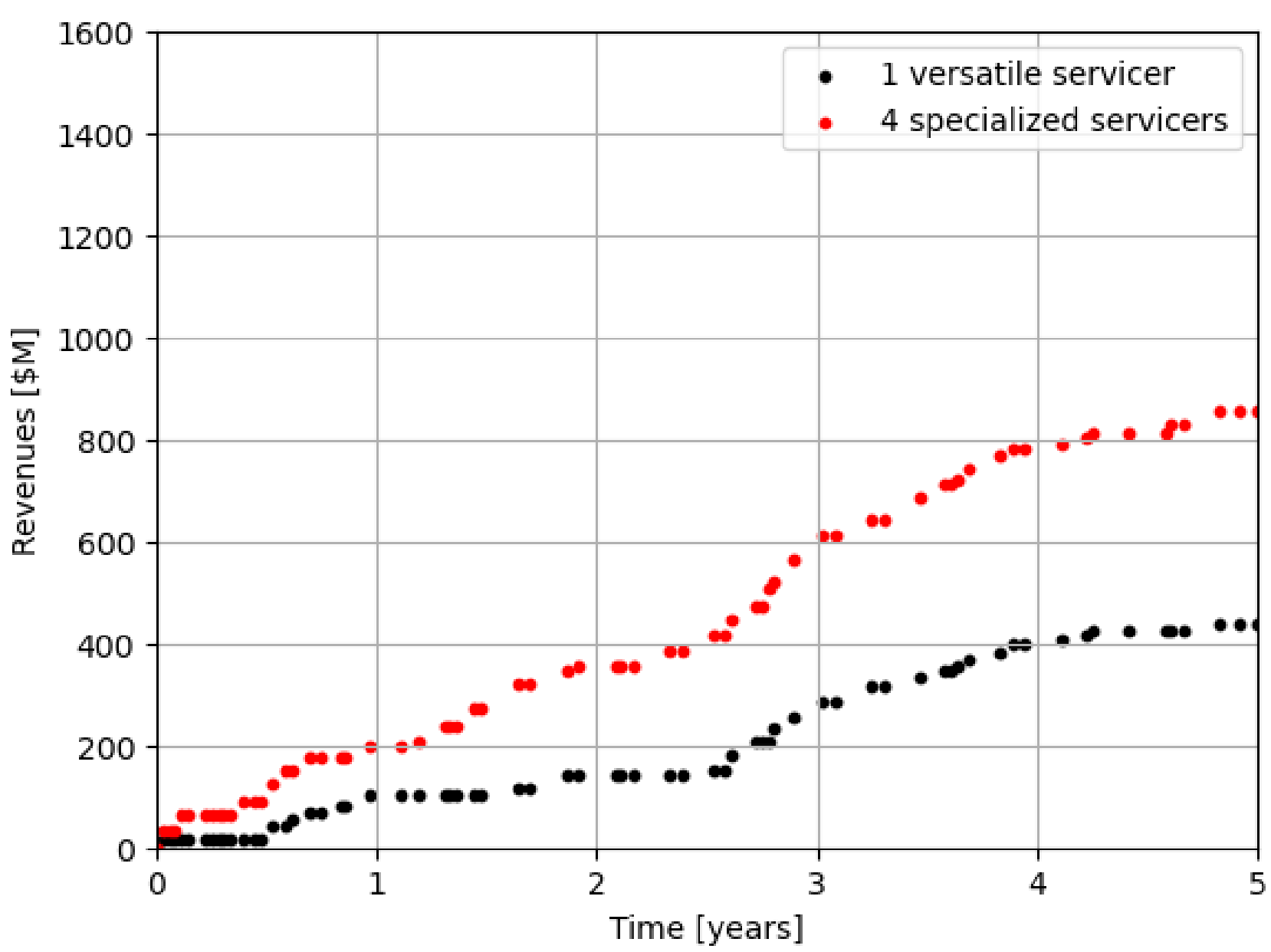

**Fig. 16 Revenues generated by the monolithic and distributed architectures for a customer base of 71 satellites.**

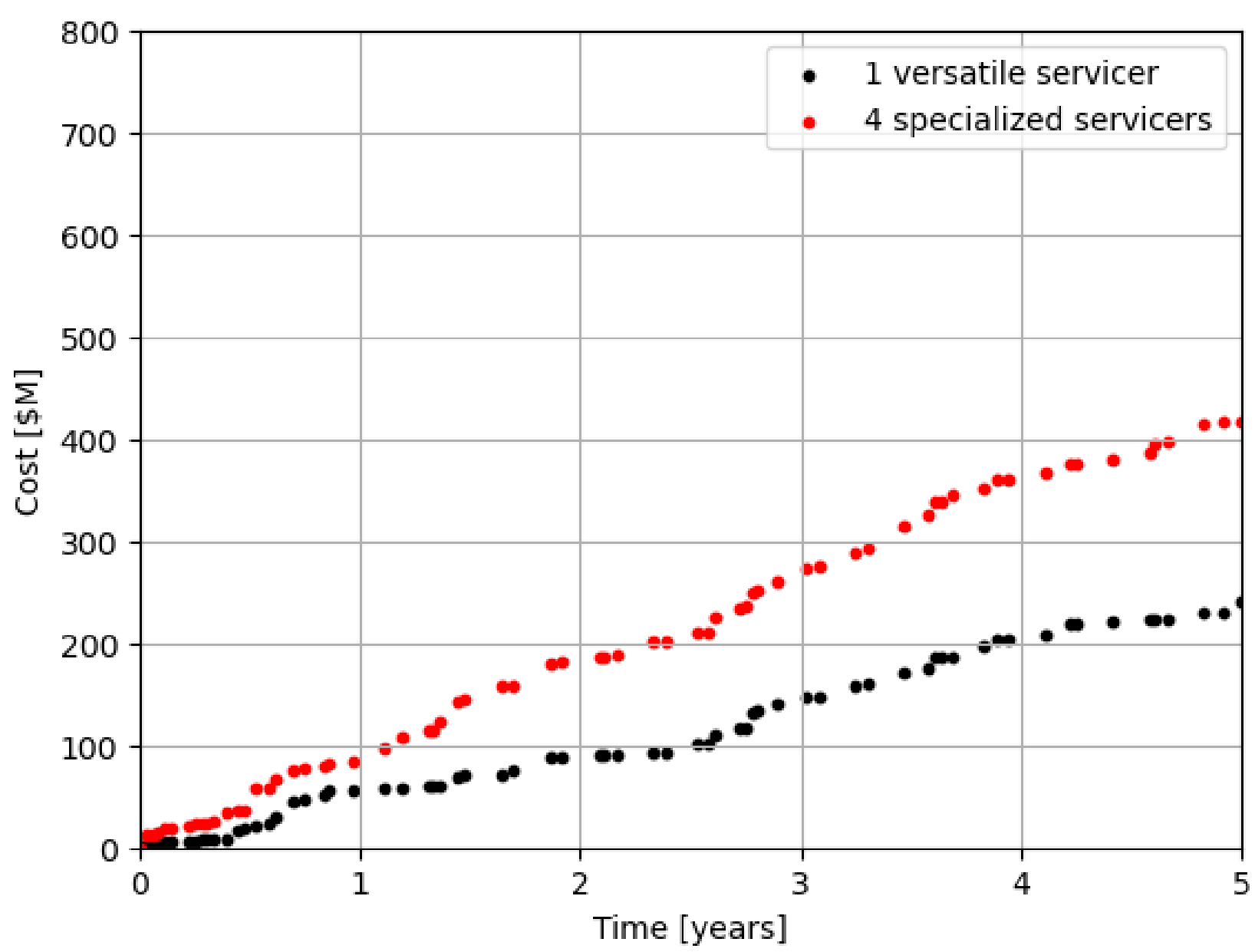

**Fig. 17 Costs incurred by the monolithic and distributed architectures for a customer base of 71 satellites.**

We finally assume that over the years the confidence in orbital servicing grows significantly. Figure 18 simulates the values of both infrastructures for a customer base of 142 satellites. Figures 19 and 20 respectively represent the revenues generated and costs incurred by both infrastructures over time. Unlike the trends observed in Fig. 15, the OOS infrastructures can pay for themselves after 2.5 years of operation for the distributed architecture and after 5 years for the monolithic architecture. However, as observed previously, the distributed architecture has a leading competitive edge over the monolithic architecture, being much more profitable.

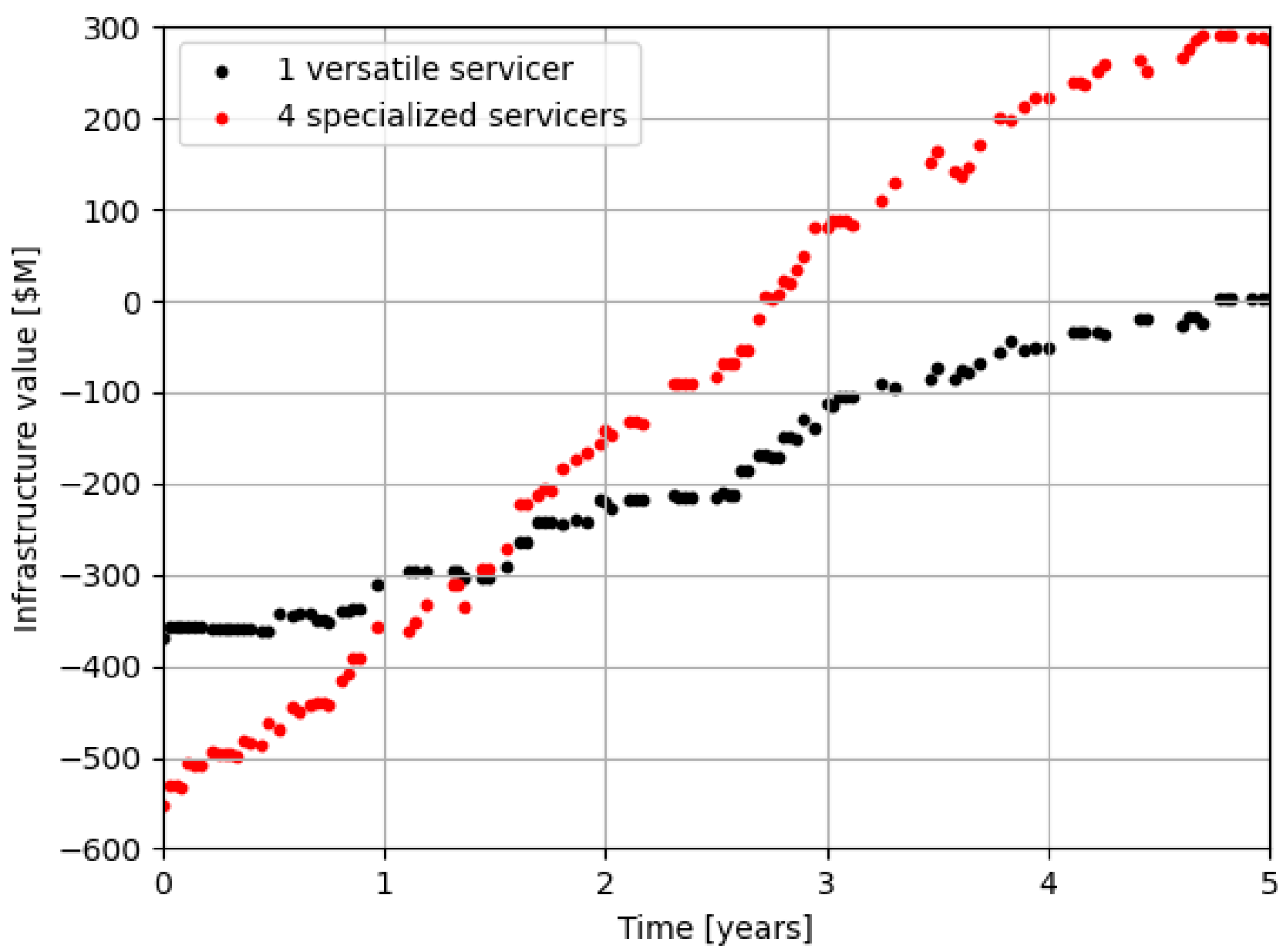

**Fig. 18 Values of the monolithic and distributed architectures for a customer base of 142 satellites.**

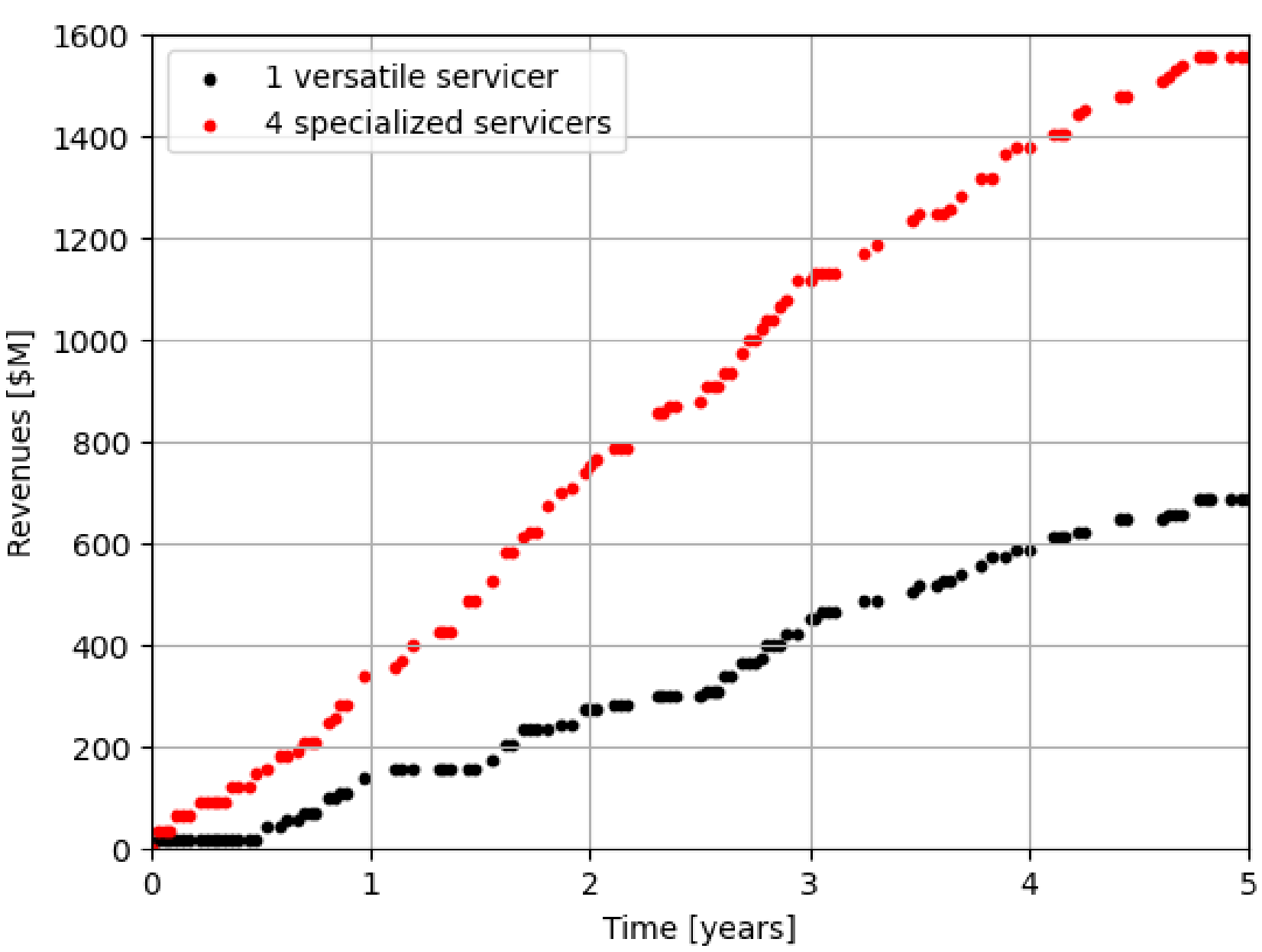

**Fig. 19 Revenues generated by the monolithic and distributed architectures for a customer base of 142 satellites.**

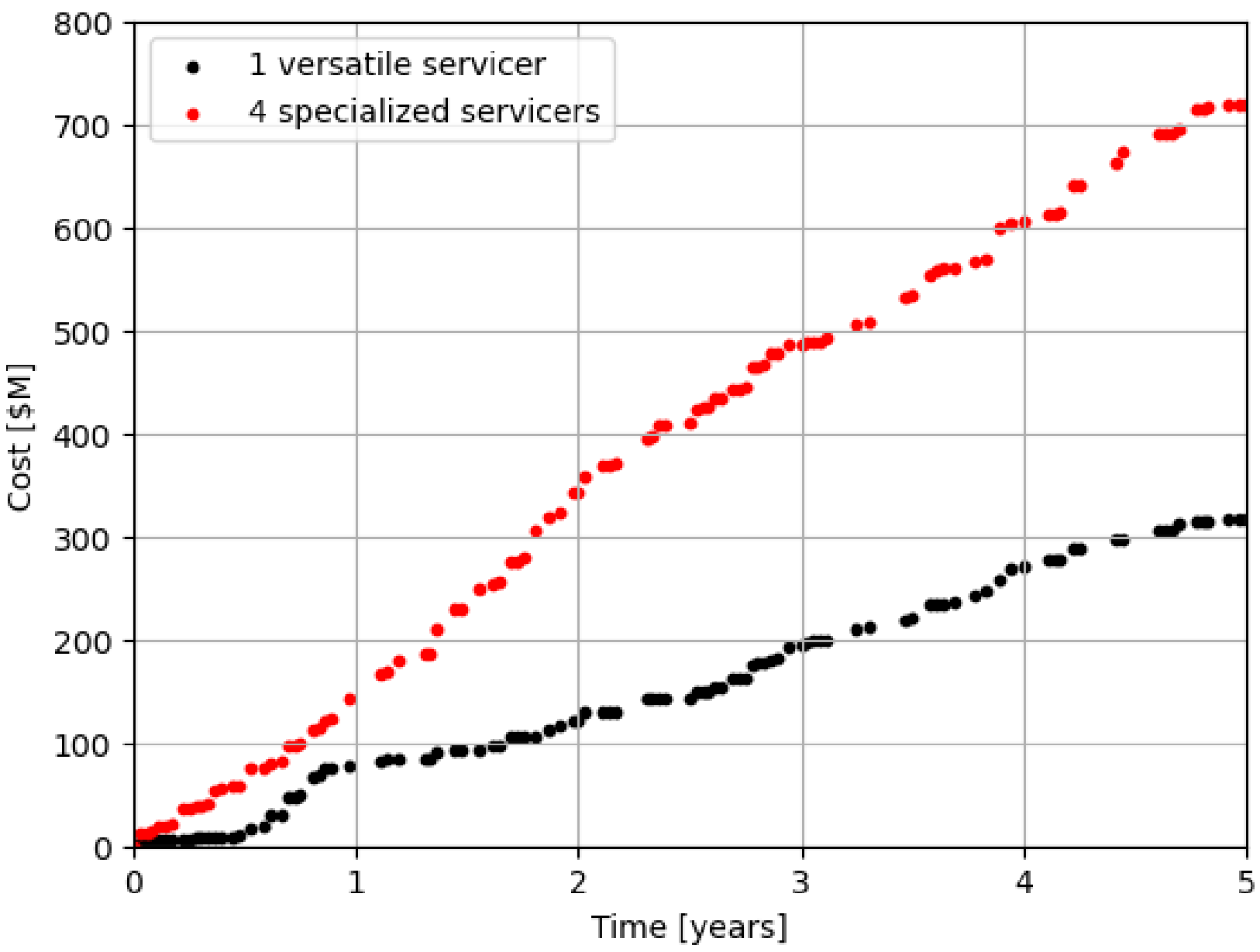

**Fig. 20 Costs incurred by the monolithic and distributed architectures for a customer base of 142 satellites**

Rather than providing definitive conclusions about the best OOS architectures given some market conditions, this section aims to demonstrate that our framework can handle the modeling and analysis of a wide spectrum of potentially complex infrastructures. A complex infrastructure was chosen to show the capability of our methodology, but the same approach can be used for a simpler, near-term infrastructure case as well; the generality is an advantage of the proposed approach.

The current framework can also easily be modified to model a dynamic evolution of the servicing demand. For example, as the confidence in OOS improves, the market size may progressively increase. This would mean that precursor OOS infrastructures could lose market shares if they reach their maximum servicing capacities. Given some assumed dynamic market forecast, sensitivity analyses could be run using our framework to gain insights into how and when to update an existing OOS infrastructure.

## V. Conclusion

This paper develops an optimization framework to enable on-orbit servicing operators to make rigorous decisions regarding the short-term operations and/or the long-term strategic planning of their fleets of servicers and orbital depots. To this end, three innovations are achieved: (1) model the operations of sustainable many-to-many OOS infrastructures as a logistics network; (2) extend the traditional space logistics MILP formulation with additional variables and constraints that model realistic OOS operations; and (3) account for the uncertainties in service demand by adapting the RH decision making approach to OOS logistics.

The proposed framework is developed based on three main assumptions: (1) the fleet of customer satellites is distributed along a shared circular orbit; (2) the servicers use impulsive-thrust technology to perform orbital maneuvers; and (3) the orbital depots, if deployed, are located on the same circular orbit as the fleet of customer satellites.

Two different case studies are explored to demonstrate the framework's value in the context of GEO servicing. The first case study shows the framework's ability in efficiently optimizing short-term operational scheduling of existing OOS infrastructures. This is done by running it over a single PH. The second case study proves the framework's value in exploring the design tradespace of OOS architectures under various OOS market conditions. This provides answers as to what features an OOS architecture should exhibit to draw the most benefits out of an assumed customer base. In this paper, three 5-year simulations are performed with three different service demand rates in order to compare two fundamentally different OOS architectures: a monolithic architecture that uses a single versatile servicer versus a distributed architecture that leverages four specialized servicers.

To conclude, the framework proposed in this paper uniquely extends the network-based space logistics technique and the RH approach, and combines them to the OOS design and operational optimization under uncertain demands. Future work will relax the assumptions made in this paper by modeling the low-thrust propulsion technology and trajectories of the servicers and by considering customer fleets distributed over different orbits.

## Appendix: Equations for Phasing Maneuvers and Rocket Dynamics

*Phasing maneuvers*

Phasing maneuver-based rendezvous are two-impulse trajectories between two different positions on the same circular orbit. Unlike Hohmann transfers, there is no waiting time before a transfer, and the travel time depends on the initial relative angle between the GEO satellite and the servicer. The lack of waiting time is what facilitates the modeling of the OOS operations as a MILP problem in this paper. The semi-major axis of the phasing orbit is calculated by minimizing the delta V under altitude and time-of-flight constraints. Figure 21 illustrates the situation when the semi-major axis of the servicer phasing trajectory is larger than the GEO orbit radius. However, the algorithm proposed here can also find optimal semi-major axes smaller than the GEO orbit radius.

The radius of the GEO orbit is denoted by $r$. Similarly, $r_{crit}$ is the radius of the "forbidden flight zone" defined around the Earth where the servicer is not allowed to fly. The initial relative angle from the servicer to the GEO satellite is $\Delta\theta_0 \in [0,2\pi)$. Similarly, the initial relative angle from the GEO satellite to the servicer is $\alpha = 2\pi - \Delta\theta_0$ and is used instead of $\Delta\theta_0$ in the equations presented below.

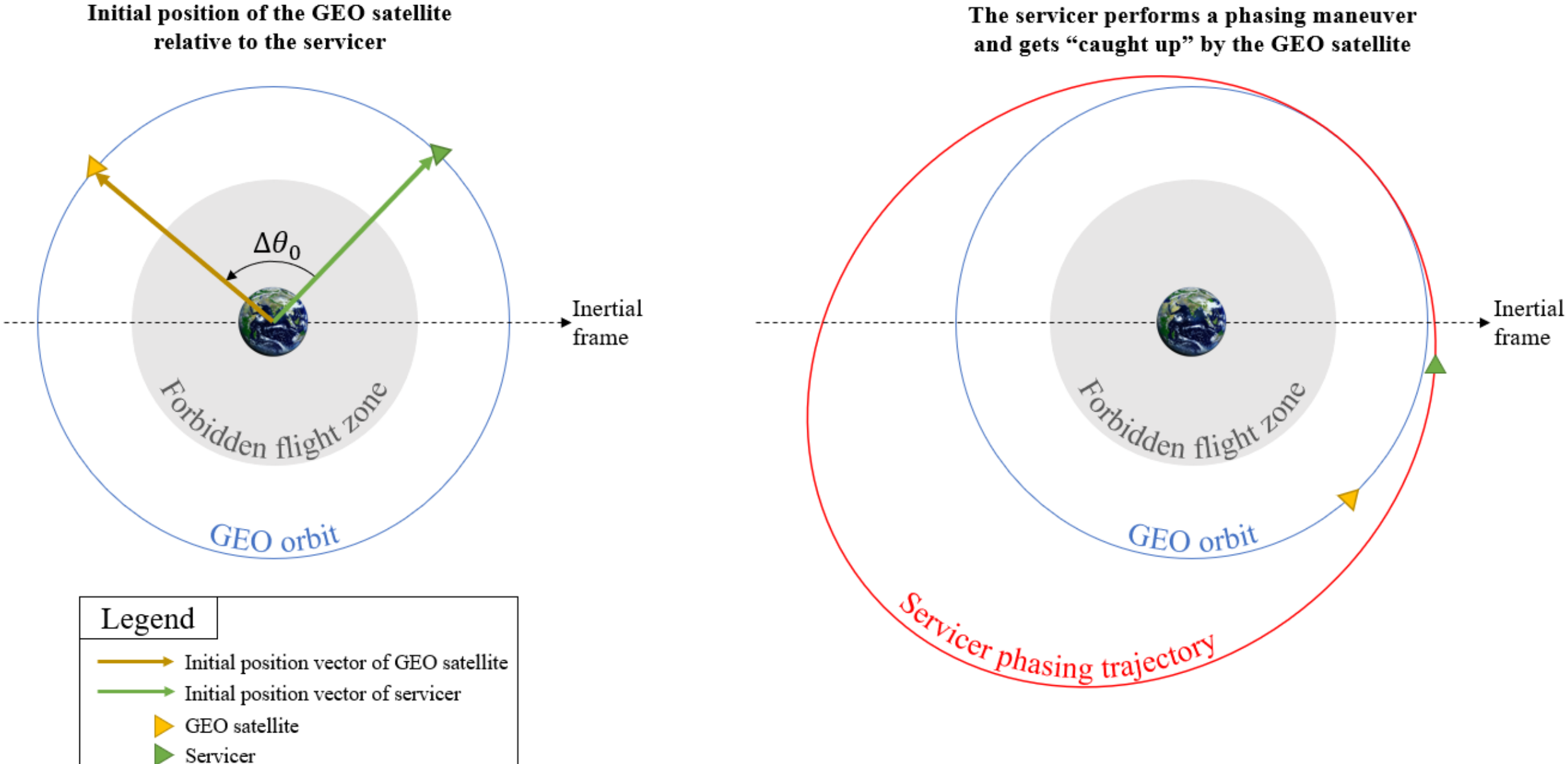

**Fig. 21 Example of rendezvous between a servicer and a GEO satellite by making use of a phasing maneuver.**

In order to define the travel time of the servicer and satellite on their respective orbit, we define two integers $k_1 \geq 1$ and $k_2 \geq 0$ corresponding to the number of complete orbits they fly before rendezvousing. From the point of view of the servicer, the travel time is

$$t_{\text{travel}} = k_1 T_1 = 2\pi k_1 \sqrt{\frac{a^3}{\mu}} \tag{21}$$

where $T_1$ and $a$ are the period and semi-major axis of the phasing orbit respectively, and $\mu$ is the gravitational parameter of the Earth. From the point of view of the satellite, the travel time is

$$t_{\text{travel}} = \frac{\alpha}{2\pi} T + k_2 T = \left[\frac{\alpha + 2\pi k_2}{2\pi}\right] T = (\alpha + 2\pi k_2) \sqrt{\frac{r^3}{\mu}} \tag{22}$$

where $T$ is the period of the orbit of the satellite. Equaling the above two equations gives an expression for $a$ as a function of $k_1$ and $k_2$:

$$a = \left[\frac{\alpha + 2\pi k_2}{2\pi k_1}\right]^{2/3} r \tag{23}$$

The altitude constraint can then be written as

$$a \geq \frac{r + r_{crit}}{2} \tag{24}$$

Solving for $k_1$ and $k_2$ is a two-step process. One can see from Eq. (22) that the travel time can be expressed as a function of only one unknown, $k_2$, from the perspective of the satellite. Thus, this equation is used to find the $k_2$ values that satisfy the time-of-flight constraint. For each tested value of $k_2$, values of $k_1$ are tested using Eq. (23) and Eq. (24). If the pair $(k_1, k_2)$ satisfies both the time-of-flight and altitude constraints, it is kept as candidate solution. Once all such feasible pairs have been found, the one with the lowest delta V is kept as solution. The delta V is then calculated using:

$$\Delta V = 2\left|\sqrt{\frac{\mu}{r}} - \sqrt{\mu\left(\frac{2}{r} - \frac{1}{a}\right)}\right| \tag{25}$$

*Rocket dynamics*

Given the delta V required for some impulsive orbital maneuver, the rocket equation is used to determine the amount of propellant needed to perform the maneuver.

$$\Delta V = g I_{sp} \ln\left(\frac{m_0}{m_0 - m_p}\right) \tag{26}$$

where $g$ is the gravitational acceleration at sea level, $I_{sp}$ is the specific impulse of the spacecraft, $m_0$ is the initial mass of the spacecraft, and $m_p$ is the propellant consumed by the spacecraft to perform the maneuver. By inverting this relationship, one can easily deduce the mass of propellant $m_p$ consumed during the maneuver.

## Acknowledgments

This work is supported by the Defense Advanced Research Project Agency Young Faculty Award D19AP00127. The content of this paper does not necessarily reflect the position or the policy of the U.S. Government, and no official endorsement should be inferred. This paper has been approved for public release; distribution is unlimited. The authors would like to thank Masafumi Isaji for his editorial review of the manuscript.